\documentclass{amsart}
\usepackage{amssymb}
\usepackage{amscd}

\newtheorem{Theorem}{Theorem}[section]
\newtheorem{Lemma}[Theorem]{Lemma}
\newtheorem{Corollary}[Theorem]{Corollary}
\newtheorem{Proposition}[Theorem]{Proposition}
\theoremstyle{definition}
\newtheorem{Definition}[Theorem]{Definition}
\newtheorem{Example}[Theorem]{Example}
\theoremstyle{remark}
\newtheorem{Remark}{Remark}

\font\sy=cmsy10

\font\ym=msbm10
\newcommand{\trace}{\hbox{\rm tr}}
\newcommand{\Aut}{\hbox{\rm Aut}}
\newcommand{\End}{\hbox{\rm End}}
\newcommand{\Hom}{\hbox{\rm Hom}}
\newcommand{\Spec}{\hbox{\rm Spec}}
\newcommand{\bA}{\text{\ym A}}
\newcommand{\bB}{\text{\ym B}}
\newcommand{\cA}{{\hbox{\sy A}}}
\newcommand{\cB}{{\hbox{\sy B}}}
\newcommand{\cC}{{\hbox{\sy C}}}

\newcommand{\cM}{{\hbox{\sy M}}}

\newcommand{\cS}{{\hbox{\sy S}}}
\newcommand{\cT}{{\hbox{\sy T}}}

\newcommand{\cV}{{\hbox{\sy V}}}

\renewcommand{\a}{\alpha}

\renewcommand{\d}{\delta}
\newcommand{\D}{\Delta}
\newcommand{\e}{\epsilon}

\newcommand{\h}{\theta}

\renewcommand{\l}{\lambda}

\newcommand{\m}{\mu}

\renewcommand{\o}{\omega}

\renewcommand{\r}{\rho}

\renewcommand{\v}{\varphi}
\newcommand{\x}{\xi}

\newcommand{\C}{{\text{\ym C}}}

\title[Tannaka Duals]
{Tannaka Duals in Semisimple Tensor Categories}
\author[Yamagami Shigeru]{Shigeru Yamagami}

\keywords{Tannaka dual, tensor category, imprimitivity bimodule, orbifold}
\subjclass{18D10, 46L37}

\begin{document}
\maketitle 
\begin{center}
Department of Mathematical Sciences\\
Ibaraki University\\
Mito, 310-8512, JAPAN\\
http://suuri.sci.ibaraki.ac.jp/\~{}yamagami/
\end{center}

\begin{abstract} 
Tannaka duals of Hopf algebras inside semisimple tensor categories are 
used to construct orbifold tensor categories, 
which are shown to include the Tannaka dual 
of the dual Hopf algebras. 
The second orbifolds are then canonically isomorphic to the initial tensor 
categories. 
\end{abstract}

\bigskip
\noindent{\bf Introduction} 

The importance of recent studies of Hopf algebras is based on its use 
as quantum symmetry, which appears more or less in connection with 
tensor categories. 
In this respect, finite group symmetry in tensor category is particularly 
interesting and provides the right place to 
take out quotients, known as the orbifold construction. 

When the relevant group is abelian, the dual group appears inside 
the orbifold tensor category in a simple way 
and we can apply the orbifold construction 
again to obtain the second orbifold tensor category, which turns out to 
recover the initial tensor category, a duality for orbifolds, in \cite{GSTC}. 

In this paper, we shall extend this kind of duality to the symmetry 
governed by Hopf algebras. 

Given a finite-dimensional semisimple Hopf algebra $A$ with its 
Tannaka dual $\cA$ realized inside a semisimple tensor category $\cT$, 
we introduce the notion of $\cA$-$\cA$ modules in $\cT$, 
which is formulated in terms of the existence of trivializing 
isomorphisms. In the group (algebra) case, this reflects the absorbing 
property of regular representations. 

The totality of our $\cA$-$\cA$ 
modules then turns out to constitute a tensor category $\cT\rtimes \cA$ 
with the unit object given by an analogue of 
the regular representation of $A$. 
The notation indicates the fact that it is a categorical analogue of 
crossed products in operator algebras. 
By the well-known crossed products vs.~fixed point algebras reciprocity, 
we may interpret $\cT\rtimes \cA$ as presenting the orbifold of $\cT$ 
by the dual Hopf algebra $A^*$. 

The orbifold tensor category $\cT\rtimes \cA$ in turn admits 
a canonical realization of the Tannaka dual $\cB$ of the dual 
Hopf algebra $A^*$, 
which allows us to take the second orbifold 
$(\cT\rtimes \cA)\rtimes \cB$ and one of our main results shows the duality 
$(\cT\rtimes \cA)\rtimes \cB \cong \cT$. 

In our previous paper \cite{GSTC}, 
we proved this for finite abelian groups 
by counting the number of simple objects 
in the second dual $(\cT\rtimes \cA)\rtimes \cB$. 
Here we shall give a more conceptual proof of duality. 
The idea has long been known in harmonic analysis of induced representations 
as imprimitivity bimodules (\cite{F}, \cite{R}). 

By forgetting the bimodule action of $\cA$ on the unit object to one-sided 
(say, right) $\cA$-action, we can make it into a right $\cB$-module $M$ with 
the property of imprimitivity, 
$M\otimes_\cB M^* \cong I$ and 
${}_\cB M^*\otimes M_\cB \cong {}_\cB I_\cB$. 

If we put $M$ into an off-diagonal corner of a suitable bicategory so that 
it connects $\cT$ and $(\cT\rtimes \cA)\rtimes \cB$, then the duality 
is obtained quite easily, though it still contains rich information on 
orbifold constructions. 

We remark here that another interesting categorical formulation of 
imprimitivity bimodules is worked out by D.~Tambara \cite{Tam}, where 
a different notion of categorical module is used to get 
an imprimitivity bimodule which relates $\cA$ and $\cB$. 

For future applications, we also investigate how the rigidity is 
inherited under the process of taking orbifolds: 
if the original tensor category $\cT$ is 
rigid and semisimple, then so is for the orbifold tensor category 
$\cT\rtimes \cA$.  

\bigskip
{\it Basic Assumptions}

We shall work with the complex number field $\C$ as a ground field, 
though any algebraically closed field of characteristic zero can be used 
equally well. 

By a tensor category, we shall mean a linear category 
with a compatible monoidal structure, which is assumed to be strict 
without losing generality by the coherence theorem. 
 
A tensor category is said to be semisimple if $\End(X) = \Hom(X,X)$ 
is a finite-dimensional semisimple algebra for any object $X$, 
which is assumed to be closed under taking subobjects and direct sums: 
To an idempotent $e$ of $\End(X)$, 
an object $eX$ (the associated subobject) 
is assigned so that $\Hom(eX,fY) = f\Hom(X,Y)e$ and  
a finite family $\{ X_j\}_{1 \leq i \leq m}$ of objects gives rise to 
an object $X_1 \oplus \dots \oplus X_m$ so that 
\[
\Hom(X_1 \oplus \dots \oplus X_m, Y_1 \oplus \dots \oplus Y_n) 
= \bigoplus_{i,j} \Hom(X_i,Y_j). 
\]

The unit object $I$ in a semisimple tensor category is assumed to be simple, 
i.e., $\End(I) = \C 1_I$, without further qualifications. 

\bigskip
\section{Bimodules in Tensor Categories} 

Let $\cT$ be a semisimple tensor category (closed under taking
subobjects and direct sums). 
By imbedding $\cT$ into $\cT\otimes \cV = \cV\otimes \cT$ with $\cV$ denoting 
the tensor category of finite-dimensional vector spaces, 
we can perform the tensor product $X\otimes V = V\otimes X$ of 
an object $X$ in $\cT$ and an object $V$ in $\cV$ so that 
\[
\Hom(X\otimes V, Y\otimes W) = \Hom(X,Y)\otimes \Hom(V,W). 
\]
Note here that the imbedding $\cT \to \cT\otimes \cV$ gives an
equivalence of tensor categories by the semisimplicity assumption on 
$\cT$. 
We also remark that, given a representative set $S$ of simple objects
in $\cT$, we have 
\[
X \cong \bigoplus_{s \in S} s\otimes \Hom(s,X)
\]
in $\cT\otimes \cV$. 

Let $A$ be a finite-dimensional semisimple Hopf algebra with 
the associated tensor category $\cA$ of finite-dimensional 
$A$-modules and consider a monoidal imbedding 
$F: \cA \to \cT$ ($F$ being a fully faithful monoidal functor). 

By a \textbf{left $\cA$-module} in $\cT$ (relative to the imbedding $F$), 
we shall mean an object $X$ in $\cT$ together with 
a natural family of isomorphisms 
$\{ \v_V: F(V)\otimes X \to X\otimes V \}$ 
(we forget the $A$-module structure of $V$, $W$ and 
regard them just vector spaces when taking the
tensor product with $X$) 
satisfying the associativity 
\[
\begin{CD}
F(V)\otimes F(W)\otimes X @>{1\otimes \v_W}>> 
F(V)\otimes X\otimes W\\
@V{m^F_{V,W}\otimes 1}VV @VV{\v_V\otimes 1}V\\ 
F(V\otimes W)\otimes X @>>{\v_{V\otimes W}}> 
X\otimes V\otimes W
\end{CD} 
\]
and the condition that 
\[
\v_\C: F(\C)\otimes X = I\otimes X \to X = X\otimes \C 
\]
is reduced to the left unit constraint $l_X$ in $\cT$. 

Let $B$ be another finite-dimensional semisimple Hopf algebra 
with $\cB$ the tensor category of $B$-modules and 
$G: \cB \to \cT$ be a monoidal imbedding. 
A \textbf{right $\cB$-module} in $\cT$ (through $G$) is, by definition, 
an object $Y$ in $\cT$ together with a natural family of isomorphisms 
$\{ \psi_W: Y\otimes G(W) \to W\otimes Y \}$ such that 
$\psi_\C = r_Y$ and 
\[
\begin{CD}
Y\otimes G(V)\otimes G(W) @>{\psi_V\otimes 1}>> 
V\otimes Y\otimes G(W)\\ 
@V{m^G_{V,W}\otimes 1}VV @VV{1\otimes \psi_W}V\\ 
Y\otimes G(V\otimes W) @>>{\psi_{V\otimes W}}> 
V\otimes W\otimes Y
\end{CD}\ .
\]

An {\bf $\cA$-$\cB$ bimodule} in $\cT$ (relative to the imbeddings $F$, $G$) 
is an object $X$ in $\cT$ together with structures of a left
$\cA$-module and a right $\cB$-module, 
\[
\v_V: F(V)\otimes X \to V\otimes X, \qquad 
\psi_W: X\otimes G(W) \to W\otimes X 
\]
such that the following diagram commutes. 
\[
\begin{CD}
F(V)\otimes X\otimes G(W) @>>> F(V)\otimes W\otimes X 
@= W\otimes F(V)\otimes X\\ 
@VVV @. @VVV\\
X\otimes V\otimes G(W) @= X\otimes G(W)\otimes V 
@>>> W\otimes X\otimes V
\end{CD}. 
\]

We shall often write ${}_\cA X_\cB$ to indicate an $\cA$-$\cB$ bimodule 
based on an object $X$ in $\cT$ when no confusion arises for 
the choice of families $\{ \v_V\}$, $\{ \psi_W \}$. 
We also use the notation 
$\x_{V,W}: F(V)\otimes X\otimes G(W) \to W\otimes X\otimes V$ to express 
the isomorphism in the above diagram, which is referred to as 
a \textbf{trivializing isomorphism} in the following. 

\begin{Example}
If $A$ is the function algebra of a finite group $H$, then 
$H$ is realized as a subset of the spectrum $\Spec(\cT)$ of $\cT$ 
through the imbedding $F$ and the functor $F$ itself is identified 
with a lift of $H \subset \Spec(\cT)$. Similarly, if $B$ is 
the function algebra of another finite group $K$, then the monoidal
imbedding $G: \cB \to \cT$ is identified with a lift of 
$K \subset \Spec(\cT)$. 

With this observation in mind, $\cA$-$\cB$ bimodules are naturally
recognized as $H$-$K$ bimodules in $\cT$. 
\end{Example}

\begin{Example}
Let $A$ be the group algebra of a finite group $G$ with 
$\cA$ the Tannaka dual of $G$. 
For notational economy, we write ${}_GV$ to express a (left) 
$G$-module with the underlying vector space $V$. 
Thus, for example, ${}_GV\otimes {}_GW$ denotes 
the tensor product $G$-module of ${}_GV$ and ${}_GW$ whereas 
${}_GV\otimes W$ means the $G$-module amplified by the vector space
$W$, with the same underlying vector space $V\otimes W$. 

Let ${}_G\C[G]$ be the left regular representation of $G$. 
Given an element $a \in G$ and a $G$-module ${}_GV$, 
define isomorphisms 
\[
\v^a_V: {}_GV\otimes {}_G\C[G] \to 
{}_G\C[G]\otimes V, \qquad 
\psi^a_V: {}_G\C[G]\otimes {}_GV \to 
V\otimes {}_G\C[G]
\]
by 
\[
\v^a_V(v\otimes g) = g\otimes ag^{-1}v, 
\qquad 
\psi^a_V(g\otimes v) = ag^{-1}v\otimes g. 
\]

Then, for any given pair $(a,b)$ of elements in $G$, 
the family $\{ \v^a_V\}$ and $\{ \psi^b_V \}$ makes ${}_G\C[G]$ 
into an $\cA$-$\cA$ bimodule in $\cA$ (relative to the trivial 
imbedding), which is denoted by ${}_\cA {R^{a,b}}_{\cA}$. 
When the left (resp.~right) action is forgotten in 
${}_\cA {R^{a,b}}_{\cA}$, the resulting left (resp.~right) $\cA$-module 
is denoted by ${}_\cA R^a$ (resp.~${R^b}_{\cA}$). 
\end{Example}

\begin{Definition}
Given Tannaka duals $\cA$, $\cB$ (of finite-dimensional semisimple
Hopf algebras) in a semisimple tensor category $\cT$ and 
$\cA$-$\cB$ bimodules ${}_\cA X_\cB$, ${}_\cA Y_\cB$ in $\cT$, 
we call a morphism $f: X \to Y$ in $\cT$ an $\cA$-$\cB$ 
{\bf intertwiner} if 
\[
\begin{CD}
F(V)\otimes X\otimes G(W) @>{1\otimes f\otimes 1}>> 
F(V)\otimes Y\otimes G(W)\\ 
@VVV @VVV\\ 
W\otimes X\otimes V @>>{1\otimes f\otimes 1}> 
W\otimes Y\otimes V
\end{CD}. 
\]

The category ${}_\cA \cT_\cB$ of $\cA$-$\cB$ bimodules in $\cT$ is 
then defined by taking $\cA$-$\cB$ intertwiners as morphisms in 
${}_\cA \cT_\cB$. 
\end{Definition}

\begin{Example}~
Let $G$ be a finite group and $\cA$ be its Tannaka dual. 
For $h \in G$, denote by $\rho(h)$ the right regular representation of 
$h$: $\rho(h): g \mapsto gh^{-1}$ for $g \in G \subset \C[G]$. 
  \begin{enumerate}
  \item 
    For $a$, $b \in G$, we have 
\[
\Hom({}_\cA R^a, {}_\cA R^b) = \C \rho(b^{-1}a) 
= \Hom({R^a}_\cA, {R^b}_\cA). 
\]
\item 
For $a'$, $b' \in G$, we have 
\[
\Hom({}_\cA {R^{a',b'}}_\cA, {}_\cA {R^{a,b}}_\cA) 
= 
\begin{cases}
\C \rho(a^{-1}a') &\text{if $a^{-1}a' = b^{-1}b'$,}\\ 
0 &\text{otherwise.}
\end{cases} 
\]
  \end{enumerate}

Recall that the underlying vector space of $R^{a,b}$ is $\C[G]$. 
\end{Example}


\section{Tensor Products} 

We shall make the totality of ${}_\cA \cT_\cB$ for various Tannaka
duals $\cA$, $\cB$ into a bicategory. 
To this end, we first introduce the notion of $\cA$-tensor products. 
Let $X_\cA$ be a right $\cA$-module and ${}_\cA Y$ be a left
$\cA$-module in $\cT$. Given a simple $A$-module $V$ and a basis 
$\{ v_i \}$ of $V$, let $\{ v_i^* \}$ be its dual basis. 
Then the linear operator $v_{i,j} = v_i\otimes v_j^*$ in $V$ is
identified with an element of $A$. These for various $V$ form 
matrix units in the algebra $A$. 
We define ${\widehat v}_{ij} \in A^*$ by 
\[
\langle {\widehat v}_{ij}, w_{kl} \rangle 
= 
\begin{cases}
\d_{il}\d_{jk} \dim V &\text{if $V \cong W$,}\\
0 &\text{otherwise.}
\end{cases} 
\]
Clearly $\{ {\widehat v}_{ij} \}_{V,i,j}$ forms a linear basis of
$A^*$. 

We now introduce an element 
$\pi({\widehat v}_{ij}) \in \End(X\otimes Y)$ by the composition 
\[
\begin{CD}
X\otimes Y @>{1\otimes \d_{F(V)}\otimes 1}>> 
X\otimes F(V^*)\otimes F(V)\otimes Y @>>> 
V^*\otimes X\otimes Y\otimes V @>>> 
X\otimes Y
\end{CD}, 
\]
where the last morphism in the diagram is given by the pairing with 
${\widehat v}_{ij}$: if the composite of the first two morphisms is 
expressed as 
\[
\sum_{i,j} v_i^*\otimes t_{ij}\otimes v_j 
\]
with $t_{ij} \in \End(X\otimes Y)$, then we set 
$\pi({\widehat v}_{ij}) = \dim(V) t_{ij}$ 
or, equivalently, 
the composite 
$X\otimes Y \to X\otimes F(V^*)\otimes F(V)\otimes Y \to 
V^*\otimes X\otimes Y\otimes V$ has the expression 
\[
\sum_{i,j} (\dim V)^{-1} v_i^*\otimes \pi({\widehat v}_{ij})\otimes
v_j, 
\]
which is an element in 
\[ 
\Hom(X\otimes Y, V^*\otimes X\otimes Y\otimes V) 
= V^*\otimes \End(X\otimes Y)\otimes V. 
\]

It is immediate to check that the map $\pi$ is consistently extended to 
the linear map of $A^*$ into $\End(X\otimes Y)$, which is again
denoted by $\pi$. 

\begin{Lemma}
Let $V$, $W$ be simple $A$-modules and $\{ v_i \}$, $\{ w_k \}$ be 
their bases. Then we have 
\[
\pi({\widehat v}_{ij}) \pi({\widehat w}_{kl}) 
= \pi({\widehat v}_{ij} {\widehat w}_{kl}). 
\]
Here the multiplication in the right hand side is the one obtained by 
dualizing the coproduct of $A$. 
\end{Lemma}

\begin{proof}
Let $
\begin{CD}
U @>{T}>> V\otimes W @>{T^*}>> U
\end{CD}
$ give an irreducible decomposition of $V\otimes W$. Then,  
for the rigidity copairing 
$\d_{V\otimes W}: \C \to W^*\otimes V^*\otimes V\otimes W$, we have 
\[
\d_{V\otimes W} = \sum_{T: U \to V\otimes W} (\overline{T}\otimes T)
\d_U, 
\]
where $\overline{T}$ is the transposed map of 
$T^*: V\otimes W \to U$. 
By the associativity and the naturality of $\cA$-actions, we see that 
the composite morphism  
\[
X\otimes Y \to X\otimes F(W^*)\otimes F(V^*)\otimes F(V)\otimes F(W) 
\otimes Y \to W^*\otimes V^*\otimes X\otimes Y 
\otimes V\otimes W
\]
is equal to 
\[
\sum_{T} 
\left( 
  \begin{CD}
    X\otimes Y @>>> U^*\otimes X\otimes Y\otimes U 
@>{\overline{T}\otimes 1\otimes T}>> 
W^*\otimes V^*\otimes X\otimes Y\otimes V\otimes W
  \end{CD}
\right), 
\]
where $X\otimes Y \to U^*\otimes X\otimes Y\otimes U$ is given by 
the composition 
\[
X\otimes Y \to X\otimes F(U^*)\otimes F(U)\otimes Y \to 
U^*\otimes X\otimes Y\otimes U. 
\]
If we replace this with 
\[
\sum_{a,b} (\dim U)^{-1} u_a^*\otimes \pi({\widehat u}_{ab}) \otimes 
u_b 
\]
and then compute $\pi({\widehat v}_{ij}) \pi({\widehat w}_{kl})$, 
we obtain the formula 
\begin{align*}
\pi({\widehat v}_{ij}) \pi({\widehat w}_{kl}) 
&= \sum_T \sum_{a,b} 
(\dim U)^{-1} 
\langle {\overline T} u_a^*\otimes \pi({\widehat u}_{ab})\otimes Tu_b, 
{\widehat v}_{ij}\otimes {\widehat w}_{kl} 
\rangle\\ 
&= \sum_T \sum_{a,b} \frac{d(V)d(W)}{d(U)} 
\langle {\overline T}u_a^*, v_i\otimes w_k \rangle 
\langle Tu_b, v_j^*\otimes w_l^* \rangle 
\pi({\widehat u}_{ab}). 
\end{align*}

On the other hand, the definition of multiplication in $A^*$ gives 
\[
\langle {\widehat v}_{ij} {\widehat w}_{kl}, x \rangle 
= \langle {\widehat v}_{ij}\otimes {\widehat w}_{kl}, \D(x) \rangle  
= \sum_T d(V)d(W) 
\langle v_j^*\otimes w_l^*, TxT^* (v_i\otimes w_k) \rangle 
\]
for $x \in A$. By using the obvious identity 
\[
T^*(v_i\otimes w_k) = \sum_a 
\langle u_a^*, T^*(v_i\otimes w_k) \rangle u_a, 
\]
the above expression takes the form 
\[
d(V)d(W) 
\sum_T \sum_a \langle v_j^*\otimes w_l^*, Txu_a \rangle 
\langle u_a^*, T^*(v_i\otimes w_k) \rangle 
\]
or equivalently we have another formula 
\[
{\widehat v}_{ij} {\widehat w}_{kl} 
= \sum_T \sum_{a,b} 
\frac{d(V)d(W)}{d(U)} 
\langle v_j^*\otimes w_l^*, Tu_b \rangle 
\langle u_a^*, T^*(v_i\otimes w_k) \rangle 
{\widehat u}_{ab}, 
\]
proving the assertion. 
\end{proof}

Since the trivial representation of $A$ is given by the counit $\e$, 
we see that $\pi(\e)$ is equal to the identity morphism as the
composition 
\[
X\otimes Y \to X\otimes I\otimes I\otimes Y 
\to \C\otimes X\otimes Y\otimes \C = X\otimes Y. 
\]
This, together with the previous lemma, shows that 
$\pi: A^* \to \End(X\otimes Y)$ is an algebra-homomorphism. 
Since $A^*$ is semisimple by Larson and Radford [LR2], the component of 
the trivial representation of $A^*$ gives rise to an idempotent 
$e_{\cA}$ in $\End(X\otimes Y)$. 
The associated subobject of $X\otimes Y$ is then denoted by 
$X\otimes_\cA Y$ and is referred to as 
the {\bf $\cA$-module tensor product} of $X$ and $Y$. 

\begin{Remark}~ 
\begin{enumerate}
\item 
The idempotent $e_\cA$ is realized as $\pi(e)$, where 
the idempotent $e$ in $A^*$ is given by 
the normalized invariant integral $e \in A^*$ of $A$: 
\[
\langle e, x \rangle 
= \sum_{[V]} \frac{\dim(V)}{\dim(A)} \trace(x_V), 
\qquad 
x \in A. 
\]
\item 
Since the counit for $A^*$ is given by the evaluation map at 
the unit $1_A$ of $A$, 
the idempotent $e_{\cA}$ is non-zero if and only if there exists 
a simple object $Z$ of $\cT$ such that 
\[
\{ f \in \Hom(Z,X\otimes Y); \pi(a^*)\circ f = 
a^*(1_A) f \quad\text{for $a^* \in A^*$} \} \not= \{ 0\}. 
\]
\end{enumerate}
\end{Remark}

Let $\cA$, $\cB$ and $\cC$ be Tannaka duals in the tensor category
$\cT$ and consider ${}_\cA X_\cB$, ${}_\cB Y_\cC$. 
The tensor product $X\otimes Y$ is then an $\cA$-$\cC$ module in an obvious
manner and the associativity of biactions for $X$, $Y$ gives the following. 

\begin{Lemma}
We have 
\[
\pi(B^*) \subset \End({}_\cA X\otimes Y_\cC). 
\]
\end{Lemma}

In particular, the biaction of $\cA$ and $\cC$ on $X\otimes Y$ is
reduced to the subobject $X\otimes_\cB Y$, which is denoted by 
${}_\cA X\otimes_\cB Y_\cC$ and is referred to as the 
{\bf relative tensor product} of bimodules. 
For morphisms $f: {}_\cA X_\cB \to {}_\cA {X'}_\cB$ and 
$g: {}_\cB Y_\cC \to {}_\cB {Y'}_\cC$, 
$f\otimes g \in \Hom({}_\cA X\otimes Y_\cC, {}_\cA {X'\otimes
  Y'}_\cC)$ obviously commutes with $\pi(B^*)$ and hence 
induces the morphism 
\[
f\otimes_\cB g: {}_\cA X\otimes_\cB Y_\cC \to 
{}_\cA {X'}\otimes_\cB {Y'}_\cC, 
\]
which is the relative tensor product of morphisms. 

The operation of taking relative tensor products is clearly
associative. Thus the categories of bimodules in $\cT$ constitute 
a bicategory if we can show the existence of unit objects.


\section{Unit Objects} 

Let $F: \cA \to \cT$ be a fully faithful imbedding of the Tannaka dual 
$\cA$ of a Hopf algebra $A$. Given $A$-modules $U$, $V$ and $W$, 
we use the notation 
\[
\begin{bmatrix}
U\\ V\,W
\end{bmatrix}
= \Hom(U, V\otimes W). 
\]
Choose a representative set $\{ V\}$ of irreducible $A$-modules and
set 
\[
\bA = \bigoplus_V F(V)\otimes V^*, 
\]
which is an object in $\cT$ (more precisely in $\cT\otimes \cV$). 
Given an $A$-module $U$, define an isomorphism 
$F(U)\otimes \bA \to \bA\otimes U$ by the composition 
\begin{align*}
F(U)\otimes \bA &= \bigoplus_V F(U)\otimes F(V)\otimes V^*\\ 
&\cong \bigoplus_V F(U\otimes V)\otimes V^* 
\quad\text{(by the multiplicativity of monoidal functor)}\\
&\cong \bigoplus_{V,X} F(X)\otimes 
\begin{bmatrix}
X\\ U\,V
\end{bmatrix} 
\otimes V^* \quad\text{(by the irreducible decomposition of $U\otimes V$)}\\ 
&\cong \bigoplus_{V,X} F(X)\otimes 
\begin{bmatrix}
V^*\\ X^*\,U
\end{bmatrix} 
\otimes V^* 
\quad\text{(by Frobenius transform)}\\ 
&= \bigoplus_X F(X)\otimes X^*\otimes U 
\quad\text{(by the irreducible decomposition of $X^*\otimes U$)}\\
&= \bA\otimes U. 
\end{align*}

Similarly, we define an isomorphism 
$\bA\otimes F(U) \to U\otimes \bA$ by 
\begin{align*}
\bA\otimes F(U) &= \bigoplus_V F(V)\otimes F(U)\otimes V^*\\ 
&\cong \bigoplus_{V,X} F(X)\otimes 
\begin{bmatrix}
X\\ V\,U
\end{bmatrix} 
\otimes V^*\\
&\cong \bigoplus_{V,X} F(X)\otimes 
\begin{bmatrix}
V^*\\ U\,X^*
\end{bmatrix} 
\otimes V^*\\
&= \bigoplus_X F(X)\otimes U\otimes X^*\\ 
&= U\otimes \bA. 
\end{align*}
Here in the last line, we applied the commutativity 
$F(X)\otimes U = U\otimes F(X)$. 

\begin{Lemma}
The isomorphisms defined so far make $\bA$ into an 
$\cA$-$\cA$ bimodule. 
\end{Lemma}

\begin{proof}
We just check the compatibility of left and right isomorphisms: 
Given $A$-modules $U$ and $W$, we shall prove the commutativity of the
diagram 
\[
\begin{CD}
F(U)\otimes \bA\otimes F(W) @>>> F(U)\otimes W\otimes \bA 
@= W\otimes F(U)\otimes \bA\\
@VVV @. @VVV\\ 
\bA\otimes U\otimes F(W) @= \bA\otimes F(W)\otimes U 
@>>> W\otimes \bA\otimes U
\end{CD}. 
\]
By the associativity of the monoidal functor $F$ 
\[
\begin{CD}
F(U)\otimes F(V)\otimes F(W) @>>> F(U)\otimes F(V\otimes W)\\ 
@VVV @VVV\\
F(U\otimes V)\otimes F(W) @>>> F(U\otimes V\otimes W)
\end{CD}, 
\]
the problem is reduced to the equality of compositions 
\begin{gather*}
\bigoplus_{V,X} F(X)\otimes 
\begin{bmatrix}
X\\ U\,V\,W
\end{bmatrix}
\otimes V^* \to 
\bigoplus_{V,X,Y} F(X)\otimes 
\begin{bmatrix}
X\\ U\,Y
\end{bmatrix} 
\otimes 
\begin{bmatrix}
Y\\ V\,W
\end{bmatrix} 
\otimes V^* \to 
\bigoplus_{V,X} F(X)\otimes 
\begin{bmatrix} 
V^*\\ WX^*U
\end{bmatrix}
\otimes V^*,\\
\bigoplus_{X,V} F(X)\otimes  
\begin{bmatrix}
X\\ U\,V\,W
\end{bmatrix}
\otimes V^* \to 
\bigoplus_{V,X,Y} F(X)\otimes 
\begin{bmatrix}
X\\ Y\,W
\end{bmatrix} 
\otimes 
\begin{bmatrix}
Y\\ U\,V
\end{bmatrix} 
\otimes V^* \to 
\bigoplus_{V,X} F(X)\otimes 
\begin{bmatrix} 
V^*\\ 
WX^*U
\end{bmatrix}
\otimes V^*. 
\end{gather*}

By an easy manipulation of transposed morphisms 
(no spherical normalization is needed here for rigidity), we see
that these are the ones associated to the following composite Frobenius
transforms 
\begin{gather*}
  \begin{bmatrix}
  X\\ U\,V\,W
  \end{bmatrix} 
\to 
\begin{bmatrix}
W^*\\ X^*\,U\,V
\end{bmatrix} 
\to 
\begin{bmatrix}
V^*\\ W\,X^*\,U
\end{bmatrix},\\ 
  \begin{bmatrix}
  X\\ U\,V\,W
  \end{bmatrix} 
\to 
\begin{bmatrix}
U^*\\ V\,W\,X^*
\end{bmatrix} 
\to 
\begin{bmatrix}
V^*\\ W\,X^*\,U
\end{bmatrix}. 
\end{gather*}

Now the coincidence of these is further reduced to the equality of 
left and right transposed morphisms, which is a consequence of 
the involutiveness of antipodes for finite-dimensional semisimple 
Hopf algebras ([LR1]). 

Given a vector 
\[
f\otimes g \in 
\begin{bmatrix}
X\\ U\,Y
\end{bmatrix} 
\otimes 
\begin{bmatrix}
Y\\ V\,W
\end{bmatrix} 
\]
in the middle vector space, we need to identify the map 
\[
\begin{bmatrix}
X\\ U\,V\,W
\end{bmatrix} 
\ni (1\otimes g)f \mapsto 
(1\otimes{\widetilde f}){\widetilde g} \in 
\begin{bmatrix}
V^*\\ W\,X^*\,U
\end{bmatrix}, 
\]
where 
\[
{\widetilde f} \in 
\begin{bmatrix}
Y^*\\ X^*\,U
\end{bmatrix}, 
\quad 
{\widetilde g} \in 
\begin{bmatrix}
V^*\\ W\,Y^*
\end{bmatrix} 
\]
are Frobenius transforms of $f$ and $g$ respectively. 
Now Fig.~\ref{biaction} shows that the morphism 
$(1\otimes {\widetilde f}){\widetilde g}$ is obatined by applying 
Frobenius transforms to $(1\otimes g)f$ repeatedly. 
\end{proof}

\begin{figure}[htbp]
\vspace*{0.2cm}
\hspace*{-1.5cm}
\input{biaction.tpc}
\caption[]{}
\label{biaction}
\end{figure}

\begin{Remark} 
We have the following gauge ambiguity for 
the choice of trivializing isomorphisms: 
Given an invertible element $\h \in \End(\bA)$, we can perturb 
the trivialization isomorphisms by the commutativity of the diagram 
\[
\begin{CD}
F(U)\otimes \bA\otimes F(W) @>{\a_{U,W}}>> W\otimes \bA\otimes U\\ 
@A{1\otimes \h\otimes 1}AA @AA{1\otimes \h\otimes 1}A\\ 
F(U)\otimes \bA\otimes F(W) @>>{\a^\h_{U,W}}> 
W\otimes \bA\otimes U
\end{CD}\ . 
\]
Note that, $\bA$ being isomorphic to $\bigoplus_V F(V)\otimes V^*$ 
as an object in $\cT$, we have the identification 
$\Aut(\bA) = \prod_V \text{GL}(V^*)$. 

When $\cT$ is a C*-tensor category and $A$ is a C*-Hopf algebra, 
with the choice of $\h$ defined by the family 
$\{ \sqrt{d(V)} 1_{V^*} \}_V$, 
the isomorphism $\a^\h_{U,W}$ becomes a unitary. 
In fact, the unperturbed isomorphism are locally given by 
\[
\begin{bmatrix}
X\\ V\,U 
\end{bmatrix} 
\otimes V^* 
\ni T\otimes v^* \mapsto {\widetilde T}v^* 
\in U\otimes X^* 
\]
with their norms (the inner products being associated to 
operator norms) by 
\[
\| T\otimes v^*\|^2 = \frac{1}{d(X)} \langle T^*T\rangle 
(v^*|v^*), 
\qquad 
\| {\widetilde T} v^* \|^2 = \frac{1}{d(V)} 
\langle T^*T \rangle (v^*|v^*). 
\]
\end{Remark}


\section{Unit Constraints} 

Given a left $\cA$-module $X$ in $\cT$, we now introduce 
a morphism $\l: \bA\otimes X \to X$ by the composition 
\[
\bigoplus_V F(V)\otimes X \otimes V^* \to 
\bigoplus_V X\otimes V\otimes V^* \to X, 
\]
where the last morphism is the one associated to the pairing map 
\[
\bigoplus_V V\otimes V^* \ni v\otimes v^* 
\mapsto \langle v, v^* \rangle \in \C. 
\]

\begin{Lemma}
We have 
\[
\l\circ \pi(a^*) = a^*(1) \l: \bA\otimes X \to X 
\]
for $a^* \in A^*$. 

Moreover, $\lambda$ is $\cA$-linear: 
\[
\begin{CD}
F(U)\otimes \bA\otimes X @>{1\otimes \l}>> F(U)\otimes X @>>> X\otimes U\\ 
@VVV @VVV @|\\ 
\bA\otimes U\otimes X @= \bA\otimes X\otimes U @>>{\l\otimes 1}> X\otimes U
\end{CD}. 
\] 
\end{Lemma}

\begin{proof}
Let $a^* = {\widetilde w}_{kl}$ be an element associated to a simple 
$A$-module $W$. 
Then the composition $\l\circ \pi({\widetilde w}_{kl})$ is given by 
\begin{align*}
\bigoplus_V F(V)\otimes V^*\otimes X 
&\to \bigoplus_V F(V\otimes W^*)\otimes F(W)\otimes X\otimes V^*\\ 
&\to 
\bigoplus_{U,V} F(U)\otimes 
\begin{bmatrix}
U\\ V\,W^*
\end{bmatrix} 
\otimes V^*\otimes X\otimes W\\ 
&\to \bigoplus_U F(U)\otimes W^*\otimes U^*\otimes W\otimes X\\ 
&\stackrel{{\widehat w}_{kl}}{\longrightarrow} 
\bigoplus_U F(U)\otimes U^*\otimes X\\ 
&\stackrel{\l}{\to} X, 
\end{align*}
which is, by the naturality of $F(\cdot)\otimes X \to
X\otimes(\cdot)$, equal to the composition 
\begin{align*}
\bigoplus_V F(V)\otimes V^*\otimes X 
&\to \bigoplus_V X\otimes V\otimes V^*\\ 
&\stackrel{1\otimes \d_W\otimes 1}{\longrightarrow} 
\bigoplus_V X\otimes V\otimes W^*\otimes W\otimes V^*\\ 
&\to \bigoplus_{U,V} X\otimes U\otimes W\otimes V^*\otimes 
\begin{bmatrix}
U\\ V\,W^* 
\end{bmatrix}\\ 
&\to \bigoplus_{U,V} 
X\otimes U\otimes W\otimes V^*\otimes 
\begin{bmatrix}
V^*\\ W^*\,U^* 
\end{bmatrix}\\ 
&\to \bigoplus_U X\otimes U\otimes W\otimes W^*\otimes U^*\\ 
&\stackrel{{\widehat w}_{kl}}{\longrightarrow} 
\bigoplus_U X\otimes U\otimes U^*\\ 
&\stackrel{\l}{\to} X. 
\end{align*}

We now compute how the operation works on vector spaces: 
\begin{align*}
v\otimes v^* 
&\mapsto \sum_m v\otimes w^*_m\otimes w_m\otimes v^*\\ 
&\mapsto \sum_{m,T,i} 
\langle (Tu_i)^*, v\otimes w_m^* \rangle 
Tu_i\otimes w_m\otimes v^*\\ 
&\mapsto 
\sum \langle (Tu_i)^*,v\otimes w_m^* \rangle 
u_i\otimes w_m\otimes {\widetilde T} v^*\\ 
&\mapsto d(W) \sum_{T,i}  
\langle (Tu_i)^*,v\otimes w_l^* \rangle 
\langle u_i\otimes w_k, {\widetilde T} v^* \rangle\\ 
&= d(W) \sum 
\langle u_i^*, T^*(v\otimes w_l^*) \rangle 
\langle u_i\otimes w_k, {\widetilde T} v^* \rangle\\ 
&= d(W) \sum_T 
\langle T^*(v\otimes w_l^*)\otimes w_k, {\widetilde T} v^* \rangle. 
\end{align*}
Here the families $\{ T: U \to V\otimes W^* \}_T$,  
$\{ T^*: V\otimes W^* \to U \}_T$ are chosen so that 
$S^*T = \d_{S,T} 1_U$ and set ${\overline T} = {}^tT^*$. 
Note that, if we denote by $\{ u_i^* \}$ the dual basis of 
$\{ u_i \}_i$, then the family $\{ {\overline T}u_i^* \}$ is 
the dual basis of the basis $\{ Tu_i\}_{T,i}$ 
of $V\otimes W^*$. 

By the relation 
\[
\sum_T {}^t{\widetilde T}(T^*\otimes 1) 
= \sum_T (1_V\otimes \e_{W^*})(TT^*\otimes 1_W) 
= 1_V\otimes \e_{W^*}, 
\]
the above operation on vector spaces ends up with 
\[
d(W) \langle v, v^* \rangle \e_{W^*}(w_l^*\otimes w_k) 
= d(W) \d_{kl} \langle v,v^* \rangle 
= {\widetilde w}_{kl}(1) \langle v, v^* \rangle. 
\]
Since the morphism $\lambda$ is associated to the pairing 
\[
v\otimes v^* \mapsto \langle v, v^*\rangle, 
\]
the above formula gives the result. 

To see the $\cA$-linearity, we again use the functoriality of 
trivializing morphisms and the problem is reduced to check the
commutativity 
\[
\begin{CD}
\bigoplus_V U\otimes V\otimes V^* @>>> U\\ 
@VVV @AAA\\ 
\bigoplus_{V,W} W\otimes 
\begin{bmatrix}
W\\ U\,V
\end{bmatrix} 
\otimes V^* 
@>>> 
\bigoplus_W W\otimes W^*\otimes U
\end{CD}, 
\]
i.e., $(1\otimes \e_V)(T\otimes 1_{V^*}) = 
(\e_W\otimes 1)(1_W\otimes {\widetilde T})$, which is an immediate
consequence of hook identities. 
\end{proof}

By the covariance just checked, the morphism $\l: \bA\otimes X \to X$ 
can be interpreted as defining 
${}_\cA \bA\otimes_\cA X \to {}_\cA X$, which is 
denoted by $l_X$. 

Conversely, consider the morphism $\m: X \to \bA\otimes X$ defined by 
\[
X \to \bigoplus_V X\otimes V\otimes V^* \to 
\bigoplus_V F(V)\otimes X\otimes V^* = \bA\otimes X, 
\]
where the first morphism is associated to the copairing 
\[
\bigoplus_V \m_V\sum_i v_i\otimes v_i^*
\]
and the weight $\{ \m_V\}$ will be specified soon after. 

Now the composition $\pi({\widetilde w}_{kl})\circ \m$ is given by 
\begin{align*}
X &\to \bigoplus_V X\otimes V\otimes V^*\\ 
&\stackrel{\d_W}{\longrightarrow} 
\bigoplus_V X\otimes V\otimes W^*\otimes W\otimes V^*\\
&\to \bigoplus_{U,V} X\otimes U\otimes 
\begin{bmatrix}
U\\ V\,W^*
\end{bmatrix} 
\otimes W\otimes V^*\\ 
&\to \bigoplus_{U,V} X\otimes U\otimes 
\begin{bmatrix}
V^*\\ W^*\,U^* 
\end{bmatrix} 
\otimes W\otimes V^*\\ 
&\to \bigoplus_U X\otimes U\otimes W\otimes W^*\otimes U^*\\ 
&\stackrel{{\widehat w}_{kl}}{\longrightarrow} 
\bigoplus_U X\otimes U\otimes U^*\\ 
&\to \bigoplus_U F(U)\otimes X\otimes U^*,  
\end{align*}
which we expect to be equal to $d(W)\d_{kl} \m$. 

To see this, we work with operations on vector spaces: 
\begin{align*}
\sum_{V,i} \m_V v_i\otimes v_i^* 
&\mapsto 
\sum_{V,i,j} \m_V v_i\otimes w_j^*\otimes w_j\otimes v_i^*\\ 
&= \sum_{V,i,j} \sum_{U,T,a} \m_V 
\langle (Tu_a)^*, v_i\otimes w_j^* \rangle 
Tu_a\otimes w_j\otimes v_i^*\\ 
&\mapsto \sum \m_V 
\langle (Tu_a)^*, v_i\otimes w_j^* \rangle 
u_a\otimes w_j\otimes {\widetilde T} v_i^*\\ 
&= d(W) \sum_{V,i} \sum_{U,T} \sum_{a,b} 
\m_V \langle (Tu_a)^*, v_i\otimes w_l^* \rangle 
\langle u_b\otimes w_k, {\widetilde T} v_i^* \rangle 
u_a\otimes u_b^*\\ 
&= d(W) \sum_{U,V,T,b} \m_V 
T^*
\left(
{}^t{\widetilde T}(u_b\otimes w_k)\otimes w_l^*
\right) 
\otimes u_b^*. 
\end{align*}
If we set $S = {}^t{\widetilde T}: U\otimes W \to V$ and 
let $S^*: V \to U\otimes W$ be the Frobenius transform of 
$T^*: V\otimes W^* \to U$, then the last expression takes the form 
\[
d(W) \sum_{U,V,S,b} 
\m_V (1\otimes \e_W)(S^*S(u_b\otimes w_k)\otimes w_l^*)\otimes u_b^*. 
\]
Applying the formula 
\[
\sum_{V,S} d(V) S^*S = d(U) 1_{U\otimes W} 
\]
for the choice $\m_V = d(V)$, 
the above summation is further reduced to 
\[
d(W) \sum_{U,b} (1\otimes \e_W)(u_b\otimes w_k\otimes w_l^*)\otimes
u_b^* 
= d(W) \d_{kl} \sum_{U,b} d(U) u_b\otimes u_b^*. 
\]

Thus, with the choice $\m_V = d(V)$, we have 
\[
\pi(a^*)\circ \m = a^*(1) \m 
\]
for $a^* \in A^*$. 

\begin{Lemma}
We now claim that 
\[
\l\circ \m = 
\left( \sum_V d(V)^2 \right) 1_X, 
\quad 
\m\circ \l = (\dim A) e_\cA = \sum_V \sum_i \pi({\widehat v}_{ii}). 
\]
\end{Lemma}

\begin{proof} 
The first relation is obvious from definitions. 

On the tensor product $\bA\otimes X$, 
the morphism $\pi({\widehat w}_{ll})$ is given by 
\begin{align*}
\bigoplus_V X\otimes V\otimes V^* 
&\to \bigoplus_V X\otimes V\otimes W^*\otimes W\otimes V^*\\ 
&\to \bigoplus_{U,V} 
U\otimes 
\begin{bmatrix}
U\\ V\,W^* 
\end{bmatrix} 
\otimes W\otimes V^*\\ 
&\to \bigoplus_{U,V} U\otimes 
\begin{bmatrix}
V^*\\ W^*\,U^* 
\end{bmatrix} 
\otimes W\otimes V^*\\ 
&\to \bigoplus_U U\otimes W\otimes W^*\otimes U^*\\ 
&\stackrel{{\widehat w}_{ll}}{\longrightarrow} 
\bigoplus_U U\otimes U^*. 
\end{align*}
According to this sequence of morphisms, we compute 
$(\dim A) e_\cA$ as follows: 
\begin{align*}
v\otimes v^* 
&\mapsto \sum_{W,k} v\otimes w_k^*\otimes w_k\otimes v^*\\ 
&\mapsto \sum 
\langle (Tu_a)^*, v\otimes w_k^* \rangle 
Tu_a\otimes w_k\otimes v^*\\ 
&\mapsto \sum 
\langle (Tu_a)^*, v\otimes w_k^* \rangle 
u_a\otimes w_k\otimes {\widetilde T}v^*\\ 
&= \sum d(W) 
\langle (Tu_a)^*, v\otimes w_l^* \rangle 
\langle u_b\otimes w_l, {\widetilde T}v^* \rangle 
u_a\otimes u_b^*\\ 
&= \sum d(W) \langle u_b\otimes w_l, {\widetilde T}v^* \rangle 
T^*(v\otimes w_l^*)\otimes u_b^*\\ 
&= \sum d(W) \langle w_l\otimes v^*, Tu_b \rangle 
T^*(v\otimes w_l^*)\otimes u_b^*\\ 
&= \sum d(W) T^*(v\otimes w_l^*)\otimes {}^tT(w_l\otimes v^*)\\ 
&= \sum d(W) (T^*\otimes {}^tT)(1\otimes \d_W\otimes 1)(v\otimes v^*). 
\end{align*} 

Now, letting $S: V^*\otimes U \to W^*$ and 
$S^*: W^* \to V^*\otimes U$ be Frobenius transforms of 
$T$ and $T^*$ respectively, we have 
\begin{align*}
\sum_{W, T} d(W) (T^*\otimes {}^tT)(1_V\otimes \d_W\otimes 1_{V^*}) 
&= \sum_{W,S} d(W) (\e_V\otimes 1_{UU^*}) 
(1_V\otimes S^*S\otimes 1_{U^*}) 
(1_{VV^*}\otimes \d_{U^*})\\ 
&= d(U) (\e_V\otimes \d_{U^*})  
\end{align*}
because of 
\[
\sum_{W,S} d(W) S^*S = d(U) 1_{V^*\otimes U}. 
\]

Thus we have 
\[ 
\sum d(W) T^*(v\otimes w_l^*)\otimes {}^tT(w_l\otimes v^*) 
= \sum d(U)\e_V(v\otimes v^*) \e_{U^*}, 
\]
which gives rise to the morphism $\m\circ \l$. 
\end{proof}

By symmetry, we may expect for the right unit constraint as well. 
Explicit computations are as follows: Define a morphism 
$\rho: X\otimes \bA \to X$ by the composition 
\[
\bigoplus_V X\otimes F(V)\otimes V^* 
\to \bigoplus_V V\otimes X\otimes V^* 
= \bigoplus_V X\otimes V\otimes V^* 
\to X, 
\]
where the last evaluation is specified by 
$v\otimes v^* \mapsto \langle v, v^* \rangle$. 
The inner morphism $\pi({\widehat w}_{kl})$ is then given by 
\begin{align*}
\bigoplus_V X\otimes F(V)\otimes V^* 
&\to \bigoplus_V X\otimes F(W^*)\otimes F(W)\otimes F(V)\otimes V^*\\ 
&\to \bigoplus_{U,V} 
X\otimes W^*\otimes F(U)\otimes 
\begin{bmatrix}
U\\ W\,V 
\end{bmatrix} 
\otimes V^*\\ 
&\to \bigoplus_{U,V} X\otimes W^*\otimes F(U)\otimes U^*\otimes W\\ 
&\to X\otimes F(U)\otimes U^*\\ 
&= X\otimes \cA. 
\end{align*}
By trivializing the functor $F$, the composition of 
$\pi({\widetilde w}_{kl})$ with the morphism $X\otimes \cA \to X$ is 
associated to the composition 
\begin{align*}
\bigoplus_V V\otimes V^*\otimes X 
&\to \bigoplus_V W^*\otimes W\otimes V\otimes V^*\otimes X\\ 
&\to \bigoplus_{U,V} 
W^*\otimes U\otimes 
\begin{bmatrix}
U\\ W\,V 
\end{bmatrix} 
\otimes V^*\otimes X\\ 
&\to \bigoplus_{U,V} 
W^*\otimes U\otimes U^*\otimes W\otimes X\\ 
&\to U\otimes U^*\otimes X\\ 
&\to X. 
\end{align*}

Now an explicit formula is obtained by working with vector spaces: 
\begin{align*}
v\otimes v^* 
&\mapsto \sum w_j^*\otimes w_j\otimes v\otimes v^*\\ 
&\mapsto \sum 
\langle (Tu_a)^*, w_j\otimes v \rangle 
w_j^*\otimes Tu_a\otimes v^*\\ 
&\mapsto \sum 
\langle (Tu_a)^*, w_j\otimes v \rangle 
w_j^*\otimes u_a\otimes {\widetilde T} v^*\\ 
&\mapsto d(W)\sum 
\langle (Tu_a)^*, w_k\otimes v \rangle 
\langle (u_b^*\otimes w_l)^*, {\widetilde T}v^* \rangle 
u_a\otimes u_b^*\\ 
&= d(W) \sum 
\langle w_l^*\otimes u_b, {\widetilde T}v^* \rangle 
T^*(w_k\otimes v)\otimes u_b^*. 
\end{align*} 

Here we shall use the identity 
\begin{align*}
\langle w_l^*\otimes u_b, \widetilde{T}v^*\rangle
&= \langle w_l^*\otimes \e_V, Tu_b\otimes v^* \rangle\\ 
&= \sum \langle v_j^*\otimes w_l^*, Tu_b \rangle 
\langle v_j, v^* \rangle\\ 
&= \langle v^*\otimes w_l^*, Tu_b \rangle 
\end{align*}
to obtain the expression 
\begin{align*}
\hspace{1cm} 
&= d(W) \sum \langle v^*\otimes w_l^*, Tu_b \rangle 
T^*(w_k\otimes v)\otimes u_b^*\\ 
&= d(W) \sum \langle {}^tT (v^*\otimes w_l^*), u_b \rangle 
T^*(w_k\otimes v)\otimes u_b^*\\ 
&= d(W) \sum T^*(w_k\otimes v)\otimes {}^tT(v^*\otimes w_l^*)\\ 
&\to d(W) \sum \e_U(T^*\otimes {}^tT) 
(w_k\otimes v\otimes v^*\otimes w_l^*)\\ 
&= d(W) \e_{WV}(TT^*\otimes 1)(w_k\otimes v\otimes v^*\otimes w_l^*)\\ 
&= d(W) \e_{WV}(w_k\otimes v\otimes v^*\otimes w_l^*)\\ 
&= d(W) \d_{kl} \langle v, v^* \rangle. 
\end{align*}

Thus $\rho\circ\pi({\widetilde w}_{kl})$ is equal to 
${\widetilde w}_{kl}(1) \rho$ and hence $\rho$ 
induces a morphism $r_X: X\otimes_\cA \bA \to X$. 

For the reverse morphism, we have 
\[
X \to \bigoplus_V X\otimes V\otimes V^* 
= \bigoplus_V V\otimes X\otimes V^* 
\to \bigoplus_V X\otimes F(V)\otimes V^*, 
\]
where the first morphism is given by 
\[
\bigoplus_V \sum_i d(V) v_i\otimes v_i^*. 
\]
Now the composition 
$X\to X\otimes \bA \to X$ is equal to 
\[
\left( 
\sum_V \dim(V)^2 
\right) 1_X
\]
whereas $X\otimes \bA \to X \to X\otimes \bA$ is given by 
\[
\left( 
\sum_V \dim(V)^2 
\right) 
e_\cA. 
\] 
Thus $r_X: X\otimes_\cA \bA \to X$ is an isomorphism of $\cA$-$\cA$ bimodules. 

\begin{Remark}
If we use the perturbed trivialization by $\a \in \Aut(\bA)$ for 
the $\cA$-$\cA$ action on $\bA$, then $\l$, $\m$ and $\rho$ are 
perturbed into $\l(\a\otimes 1)$, $(\a^{-1}\otimes 1)\mu$ and 
$\rho(1\otimes \a)$ respectively. 

In particular, if $\cT$ is a C*-tensor category, we obtain unitary 
constraints by taking $\a = \{ \sqrt{d(V)} 1_{V^*} \}_V$, i.e., 
they are associated to the pairing (copairing) 
\begin{gather*}
V\otimes V^* \ni v\otimes v^* \mapsto 
\sqrt{\d(V)} \langle v, v^* \rangle,\\ 
\sqrt{d(V)} \sum_i v_i\otimes v_i^* \in V\otimes V^*. 
\end{gather*}
\end{Remark}


\section{Triangle Identities} 

We shall now check the triangle identity for $\{ l_X, r_X \}$, i.e., 
given $\cA$-modules $X_\cA$ and ${}_\cA Y$, the idempotent 
$e_\cA \in \End(X\otimes Y)$ equalizes $\rho\otimes 1$ and $1\otimes \lambda$ as  
\[
X\otimes \bA\otimes Y\ \mathop{\longrightarrow}_{1\otimes \l}^{\r\otimes 1}\  
X\otimes Y\   
\stackrel{e_\cA}{\longrightarrow} \ X\otimes Y. 
\]

By the formula 
\[
e_\cA = \frac{1}{\dim A} 
\sum_{U,i} \pi({\widehat u}_{ii}), 
\]
we need to consider the composition of 
\[
\begin{CD}
X\otimes F(V)\otimes V^*\otimes Y @>>> V\otimes X\otimes V^*\otimes Y 
@>>> X\otimes Y\\ 
@| @. @|\\ 
X\otimes F(V)\otimes V^*\otimes Y @>>> X\otimes V^*\otimes Y\otimes V 
@>>> X\otimes Y
\end{CD} 
\]
with 
\[
\begin{CD}
XY @>{\oplus 1\otimes \d_{F(W)}\otimes 1}>> \bigoplus_W XF(W^*)F(W)Y 
@>>> \bigoplus_W W^*XYW 
@>{\sum {\widehat w}_{kk}}>> XY. 
\end{CD}
\]

By the associativity of trivialization, we are faced to compare 
\begin{equation}
  \begin{CD}
    XF(V)Y @>>> \bigoplus_W X F(V) F(W^*) F(W) Y 
@>>> \bigoplus_W VW^*XYW 
@>{\sum {\widehat w}_{kk}}>> VXY
  \end{CD}
\end{equation} 
and 
\begin{equation}
  \begin{CD}
    XF(V)Y @>>> \bigoplus_U X F(U) F(U^*) F(V) Y 
@>>> \bigoplus_U UXYU^*V 
@>{\sum {\widehat u}_{ii}}>> XYV
  \end{CD}
\end{equation} 
with the identification $V\otimes X\otimes Y = X\otimes Y\otimes V$. 

To this end, we choose the diagram 
\[
\begin{CD}
F(V) @>{1\otimes \d_{F(V)}}>> 
\bigoplus_W F(V)\otimes F(W^*)\otimes F(W)\\ 
@V{\d_{F(U)^*}\otimes 1}VV @VVV\\ 
\bigoplus_U F(U)\otimes F(U^*)\otimes F(V) @<<< 
\bigoplus_{U,W} F(U)\otimes 
\begin{bmatrix}
U\\ V\,W^* 
\end{bmatrix} 
\otimes F(W) 
\end{CD}
\]
so that it is ommutative, where the right vertical arrow is given by 
an irreducible decomposition 
$\{ 
\begin{CD}
F(U) @>{T}>> F(V)\otimes F(W^*) @>{T^*}>> F(U) 
\end{CD}
\}$ 
and the bottom line by an irreducible decomposition 
$\{ 
\begin{CD}
F(W) @>{S}>> F(U^*)\otimes F(V) @>{T^*}>> F(W) 
\end{CD}
\}$. 

The diagram is commutative if $S$ and $T$ are related 
so that 
\[
S = \frac{d(W)}{d(U)} \widetilde T
\]
with $\widetilde T$ the Frobenius transform of $T$. 
In fact, the relation ensures the identity 
\[
\sum_T (T^*\otimes S)(1_V\otimes \d_W) = \d_{U^*}\otimes 1_V. 
\]

By sandwiching the above diagram by $X\otimes \cdot \otimes Y$ and 
then applying trivialization isomorphisms, we obtain the commutative 
diagram 
\[
\begin{CD}
XF(V)Y @>>> \bigoplus_W XF(V)F(W^*)F(W)Y 
@>>> \bigoplus_{U,W} XF(U) 
\begin{bmatrix}
U\\ V\,W^*
\end{bmatrix} 
F(W)Y\\ 
@. @VVV @VVV\\ 
@. \bigoplus_W VW^*XYW @>>> 
\bigoplus_{U,W} UX 
\begin{bmatrix}
U\\ V\,W^*
\end{bmatrix} 
YW
\end{CD} 
\]
\[
\begin{CD} 
\bigoplus_{U,W} XF(U) 
\begin{bmatrix}
U\\ V\,W^*
\end{bmatrix} 
F(W)Y 
@>>> 
\bigoplus_U XF(U)F(U^*V)Y @. \\
@VVV @VVV @.\\
\bigoplus_{U,W} 
UX 
\begin{bmatrix}
U\\ V\,W^* 
\end{bmatrix} 
YW 
@>>> \bigoplus_U UXYU^*V 
@>>> XYV, 
\end{CD}
\]
where the upper route is exactly the morphism (1). 

To identify the lower route, we inspect the morphism 
\[
\begin{CD}
\bigoplus_W V W^*W @>>> 
\bigoplus_{U,W} U 
\begin{bmatrix}
U\\ V\,W^* 
\end{bmatrix} 
W 
@>>> 
\bigoplus_U UU^*V @>>> V
\end{CD}
\]
as 
\begin{align*}
v\otimes w^*\otimes w 
&\mapsto \sum \langle (Tu_i)^*, v\otimes w^* \rangle 
Tu_i\otimes w\\ 
&\mapsto \sum T^*(v\otimes w^*)\otimes Sw\\ 
&\mapsto \sum d(U) (\e_U\otimes 1) (T^*\otimes S) 
(v\otimes w^*\otimes w). 
\end{align*} 

The last summation is computed with the help of the relation 
\begin{align*}
\sum_{U,T} d(U) (\e_U\otimes 1) (T^*\otimes S) 
&= \sum d(W) (\e_U\otimes 1) (T^*\otimes {\widetilde T})\\ 
&= d(W) \sum (1\otimes \d_W) (TT^*\otimes 1_W)\\ 
&= d(W) 1_V\otimes \d_W 
\end{align*} 
to get $\langle w^*, w\rangle v$, which is equal to 
\[
\sum_{U,i} \langle {\widetilde u}_{ii}, w^*\otimes w \rangle v. 
\]
Thus the bottom route turns out to be the composition 
\[
\begin{CD}
XF(V)Y @>>> VXY @>>> \bigoplus_W VX F(W^*) F(W) Y @>>> 
\bigoplus_W VW^*XYW @>{\sum {\widehat w}_{kk}}>> VXY, 
\end{CD} 
\]
showing the equality of the morphisms (1), (2). 

As a summary, we conclude here the following.  

\begin{Proposition} Given a semisimple tensor category $\cT$, 
we have constructed the semisimple bicategory $\cM(\cT)$ of 
bimodules indexed by Tannaka duals of finite-dimensional 
semisimple Hopf algebras realized in $\cT$. 
More precisely, given a family $\{ \o_A \}$ of weights indexed by Hopf algebras 
realized inside $\cT$, the pair $(\o_A l_X, \o_B r_X)$ 
with $X = {}_\cA X_\cB$ gives unit constraints. 
\end{Proposition}

\begin{Remark} 
Given a Tannaka dual $\cA$ in $\cT$, 
it is not obvious, at first glance, how big is the tensor category 
${}_\cA \cM(\cT)_\cA$ of $\cA$-$\cA$ bimodules. 

It turn out in \S~7 to be large enough to recover the initial tensor category 
because $\cT$ is realized as the tensor category of $\cB$-$\cB$ bimodules 
in ${}_\cA \cM(\cT)_\cA$ with the Tannaka dual $\cB$ of the dual Hopf 
algebra $A^*$ being imbedded into ${}_\cA \cM(\cT)_\cA$ (see Theorem~7.5). 
\end{Remark} 

\begin{Lemma} 
Let $A$ be a finite-dimensional semisimple Hopf algebra 
with $\cA$ the tensor category of finite-dimensional $A$-modules. 
Given an imbedding $F: \cA \to \cT$ of $\cA$ into a semisimple tensor category 
$\cT$, let $\bA = \bigoplus_V F(V)\otimes V^*$ be the associated object. 
Then both of ${}_\cA \bA$ and $\bA_\cA$ are irreducible as $\cA$-modules. 
\end{Lemma}

\begin{proof}
Let 
\[
\phi = \bigoplus_V \phi_{V^*} \in 
\bigoplus_V \cB(V^*) = \End(\bA)
\]
belong to $\End({}_\cA \bA)$, i.e., 
\[
\begin{CD}
F(U)\otimes \bA @>>> \bA\otimes U\\ 
@V{1\otimes \phi}VV @VV{\phi\otimes 1}V\\ 
F(U)\otimes \bA @>>> \bA\otimes U
\end{CD}
\]
for any $U$. The commutativity is then equivalent to 
\[
\begin{CD}
\oplus_{V,W} F(W)\otimes 
\begin{bmatrix}
W\\ U\,V
\end{bmatrix} 
\otimes V^* 
@>>> \oplus_W F(W)\otimes W^*\otimes U\\ 
@VVV @VVV\\ 
\oplus_{V,W} F(W)\otimes 
\begin{bmatrix}
W\\ U\,V
\end{bmatrix} 
\otimes V^* 
@>>> \oplus_W F(W)\otimes W^*\otimes U
\end{CD}.  
\] 
Removing the $F(W)$ factor, we have 
\[
\begin{CD}
  \begin{bmatrix}
  V^*\\ 
W^*\,U
  \end{bmatrix} 
\otimes V^* @>>> W^*\otimes U\\ 
@V{1\otimes \phi}VV @VV{1\otimes \phi}V\\ 
\begin{bmatrix}
V^*\\ W^*\,U
\end{bmatrix} 
\otimes V^* 
@>>> W^*\otimes U
\end{CD}
\] 
for any $U$, $V$ and $W$, 
which means the equality 
\[ 
T\phi_{V^*} = (\phi_{W^*}\otimes 1_U)T 
\]
for any $T: V^* \to W^*\otimes U$. 

If we take $V = \C$ and $U = W$ with $T = \d_W$, then the condition is reduced to 
\[
\phi_\C \sum_k w_k^*\otimes w_k = \sum_k \phi_{W^*} w_k^*\otimes w_k,  
\] 
which is equivalent to $\phi_\C w_k^* = \phi_{W^*} w_k^*$ for any $k$, i.e., 
$\phi_{W^*} = \phi_\C 1_{W^*}$ for any $W$. 
Thus, it is proportional to the identity morphism $1_{\bA}$. 
\end{proof}

\begin{Remark}
The triangle identities are satisfied for perturbed 
$\cA$-$\cA$ actions on $\bA$ as well. 
Particularly, when $\cT$ is a C*-tensor category, 
the unitary constraints for the choice $\h = \{ \sqrt{d(V)} 1_{V^*} \}$ 
of perturbation satisfy the triangle identity and hence give rise to 
unit objects, i.e., $\cM(\cT)$ is a C*-bicategory. 
\end{Remark}

Finally we record here that, other than the perturbation for actions,
there remains somewhat trivial 
freedom for the choice of unit constraints: given a family 
$\{ \o_A \}_A$ of non-zero scalars, the unit constarints 
$l_X: {}_\cA \bA\otimes_\cA X_\cB \to {}_\cA X_\cB$ and 
$r_X: {}_\cA X\otimes_\cB \bB_\cB \to {}_\cA X_\cB$ are modified by 
multiplying $\o_A$ and $\o_B$ respectively.


\section{Rigidity} 

Let ${}_\cA X_\cB$ be an $\cA$-$\cB$ module in $\cT$ and 
suppose that $X$ admits a dual object $X^*$ with a rigidity pair 
$\e_X: X\otimes X^* \to I$, $\d_X: I \to X^*\otimes X$. 
On the image of $\cA$ in $\cT$, we have the natural choice of 
dual objects (and rigidity pairs), which enables us to define 
rigidity pairs such as $\e_{F(V)X} = \e_{F(V)}(1\otimes \e_X\otimes 1)$, 
$\d_{F(V)X} = (1\otimes \d_{F(V)}\otimes 1)\d_X$. Note here that 
the rigidity for $F(V)$ satisfies the Frobenius duality and we can freely 
use the relation such as $F(V)^{**} = F(V)$ while we should be careful 
when the object $X$ is involved 
because there is no privileged identification. 

Now, by applying the operation of taking transposed morphisms, 
we make $X^*$ into a $\cB$-$\cA$ module: 
the trivializing isomorphism 
$G(W)\otimes X^*\otimes F(V) \to V\otimes X^*\otimes W$ is defined 
to be the transposed morphism of the isomorphism 
$\phi: W^*\otimes X\otimes V^* \to F(V^*)\otimes X\otimes G(W^*)$: 
\[
(1_{VX^*W}\otimes \e_{F(V^*)XG(W^*)}) 
\phi
(\d_{W^*XV^*}\otimes 1_{G(W)X^*F(V^*)}). 
\]

\begin{Lemma}
We have the commutative diagrams 
\[
\begin{CD}
X\otimes X^*\otimes F(V) @>>> X\otimes V\otimes X^*\\ 
@V{\e\otimes 1}VV @VVV\\ 
F(V) @<<{1\otimes \e}< F(V)\otimes X\otimes X^*
\end{CD}, 
\qquad 
\begin{CD}
G(W)\otimes X^*\otimes X @>>> X^*\otimes W\otimes X\\ 
@A{1\otimes \d}AA @VVV\\ 
G(W) @>>{\d\otimes 1}> X^*\otimes X\otimes G(W)
\end{CD}. 
\]
\end{Lemma}

\begin{proof}
The composite morphism 
$X\otimes X^*\otimes F(V) \to X\otimes V\otimes X^* 
\to F(V)\otimes X\otimes X^* \to F(V)$ is given by 
\[
(1_{F(V)}\otimes \e_X \otimes \e_{F(V^*)X})
(\v_V^{-1}\otimes 1_{X^*}\otimes \v_{V^*}^{-1}) 
(1_X\otimes \d_{XV^*}), 
\]
where the hook identity is used to get the expression 
\[
(1_{F(V)}\otimes \e_{F(V^*)X}) 
(1_{F(V)}\otimes \v_{V^*}^{-1}) 
(\v_V^{-1}\otimes 1_{V^*}) 
(1_X\otimes \d_{V^*}). 
\]
Now we apply the associativity of $\v$, 
$\v_{V\otimes V^*} = 
(\v_V\otimes 1)(1\otimes \v_{V^*})$, 
to obtain 
\[
(1_{F(V)}\otimes \e_{F(V^*)X})(\d_{F(V^*)}\otimes 
1_{XX^*F(V)}) = \e_X\otimes 1_{F(V)}. 
\]
\end{proof}

\begin{Corollary}
The following diagrams commute 
\[
\begin{CD}
X\otimes X^* @>>> X\otimes V\otimes V^*\otimes X^* @>>> 
F(V)\otimes X\otimes X^*\otimes V^*\\ 
@VVV @. @VVV\\ 
X\otimes V^*\otimes V\otimes X^* 
@>>> V^*\otimes X\otimes X^*\otimes F(V) 
@>>> V^*\otimes F(V) = F(V)\otimes V^*
\end{CD},  
\]
\[
\begin{CD}
W^*\otimes G(W) = G(W)\otimes W^* @>>> G(W)\otimes X^*\otimes X\otimes W^* 
@>>> X^*\otimes W\otimes W^*\otimes X\\ 
@VVV @. @VVV\\ 
W^*\otimes X^*\otimes X\otimes G(W) @>>> X^*\otimes W^*\otimes W\otimes X 
@>>> X^*\otimes X
\end{CD}. 
\]
\end{Corollary}

Define the morphism 
\[
\e: X\otimes X^* \to \bA = \bigoplus_V F(V)\otimes V^*
\]
by the weighted summation of the above morphisms over $[V]$ with 
weight $\dim V$. Similarly we introduce the morphism 
\[
\d: \bB = \bigoplus_W G(W)\otimes W^* \to X^*\otimes X 
\]
by taking the summation on $[W]$ without weights. 

\begin{Lemma}
The morphism $\e: X\otimes X^* \to \bA$ is $\cA$-$\cA$ linear, whereas 
$\d: \bB \to X^*\otimes X$ is $\cB$-$\cB$ linear. 
\end{Lemma}

\begin{proof}
Consider the commutativity of the diagram 
\[
\begin{CD}
F(U)\otimes X\otimes X^* @>>> F(U)\otimes \bA\\ 
@VVV @VVV\\ 
X\otimes X^*\otimes U @>>> \bA\otimes U
\end{CD}. 
\]
The composite morphism 
$F(U)\otimes X\otimes X^* \to F(U)\otimes \bA \to \bA\otimes U$ is given by 
\begin{align*}
F(U)\otimes X\otimes X^*  
&\to \bigoplus_V F(U)\otimes X\otimes V\otimes V^*\otimes X^*\\ 
&\to \bigoplus_V F(U)\otimes F(V)\otimes X\otimes X^*\otimes V^*\\ 
&\to \bigoplus_V F(U)\otimes F(V)\otimes V^*\\ 
&\to \bigoplus_{V,W} F(W)\otimes 
\begin{bmatrix}
W\\ U\,V 
\end{bmatrix} 
\otimes V^*\\ 
&\to \bigoplus_W F(W)\otimes W^*\otimes U. 
\end{align*} 

By the naturality of the trivialization 
$F(\cdot)\otimes X \to X\otimes (\cdot)$, 
this composition can be described by 
\begin{align*}
F(U)\otimes X\otimes X^* 
&\to X\otimes U\otimes X^*\\ 
&\to \bigoplus_V X\otimes U\otimes V\otimes V^*\otimes X^*\\ 
&\to \bigoplus_{V,W} X\otimes W\otimes 
\begin{bmatrix}
W\\ U\,V 
\end{bmatrix} 
\otimes V^*\otimes X^*\\ 
&\to \bigoplus_W X\otimes W^*\otimes U\otimes X^*\\ 
&\to \bigoplus_W F(W)\otimes X\otimes X^*\otimes W^*\otimes U\\ 
&\to \bigoplus_W F(W)\otimes W^*\otimes U, 
\end{align*}
whence the problem is reduced to showing 
\[
\begin{CD}
U @>>> \bigoplus_V U\otimes V\otimes V^*\\ 
@VVV @VVV\\ 
\bigoplus_W W\otimes W^*\otimes U 
@<<< 
\bigoplus_{V,W} W\otimes 
\begin{bmatrix}
W\\ U\,V
\end{bmatrix}
\otimes V^*
\end{CD}. 
\]

The commutativity of this diagram is then a routine work of 
Frobenius transforms: The longer circuit is given by 
\begin{align*}
u &\mapsto \sum_{V,j} d(V) u\otimes v_j\otimes v_j^*\\ 
&\mapsto \sum_{T,W,k} d(V) 
\langle (Tw_k)^*, u\otimes v_j \rangle Tw_k\otimes v_j^*\\ 
&\mapsto \sum d(V) \langle w_k^*, T^*(u\otimes v_j) \rangle 
w_k\otimes {\widetilde T}v_j^*\\ 
&= \sum d(V) T^*(u\otimes v_j)\otimes 
{\widetilde T}v_j^*. 
\end{align*}

By replacing the summation indices $T$ and $T^*$ by their Frobenius transforms 
$S: U^*\otimes W \to V$ and $S^*: V \to U^*\otimes W$, 
we have  
\begin{align*}
\sum_{T,V} d(V) (T^*\otimes {\widetilde T}) 
(1_U\otimes \d_{V^*}) 
&= \sum_{S,V} d(V) (\e_U\otimes 1_W) 
(1_U\otimes S^*S\otimes 1_{W^*U}) 
(1_U\otimes \d_{W^*U})\\ 
&= d(W) (\e_U\otimes 1_W)(1_U\otimes \d_{W^*U})\\ 
&= d(W) \d_{W^*}\otimes 1_U,  
\end{align*}
which is used to get 
\[
\sum d(V) T^*(u\otimes v_j)\otimes {\widetilde T} v_j^* 
= \sum_{W,k} d(W) w_k\otimes w_k^*\otimes u. 
\]

A bit of care is needed for the right action: 
\[
\begin{CD}
X\otimes X^*\otimes F(U) @>>> \bA\otimes F(U)\\ 
@VVV @VVV\\ 
U\otimes X\otimes X^* @>>> U\otimes \bA
\end{CD}. 
\]
By using the previous lemma, the composite morphism 
$X\otimes X^*\otimes F(U) \to \bA\otimes F(U) \to U\otimes \bA$ is given by 
\begin{align*}
X\otimes X^*\otimes F(U) &\to 
\bigoplus_V X\otimes V^*\otimes V\otimes X^*\otimes F(U)\\ 
&\to \bigoplus_V V^*\otimes X\otimes X^*\otimes F(V)\otimes F(U)\\ 
&\to \bigoplus_V V^*\otimes F(V)\otimes F(U)\\ 
&\to \bigoplus_{V,W} V^*\otimes F(W)\otimes 
\begin{bmatrix}
W\\ V\,U
\end{bmatrix}\\ 
&\to \bigoplus_W U\otimes W^*\otimes F(W). 
\end{align*} 
By the naturality of trivialization, this is equal to 
\begin{align*}
X\otimes X^*\otimes F(U) 
&\to X\otimes U\otimes X^*\\ 
&\to \bigoplus_V X\otimes V^*\otimes V\otimes U\otimes X^*\\ 
&\to \bigoplus_{V,W} 
X\otimes V^*\otimes W\otimes 
\begin{bmatrix}
W\\ V\,U 
\end{bmatrix} 
\otimes X^*\\ 
&\to \bigoplus_W X\otimes U\otimes W^*\otimes W\otimes X^*\\ 
&\to \bigoplus_W U\otimes W^*\otimes X\otimes X^*\otimes F(W)\\ 
&\to \bigoplus_W U\otimes W^*\otimes F(W). 
\end{align*}
If we compare this with the other composite morphism 
\begin{align*}
X\otimes U\otimes X^* 
&\to \bigoplus_W X\otimes U\otimes W^*\otimes W\otimes X^*\\ 
&\to \bigoplus_W U\otimes W^*\otimes X\otimes X^*\otimes F(W)\\ 
&\to \bigoplus_W U\otimes W^*\otimes F(W), 
\end{align*}
then the problem is reduced to the commutativity of 
\[
\begin{CD}
U @>>> \bigoplus_V V^*\otimes V\otimes U\\ 
@VVV @VVV\\ 
\bigoplus_W U\otimes W^*\otimes W 
@<<< 
\bigoplus_{V,W} V^*\otimes W\otimes 
\begin{bmatrix}
W\\ V\,U
\end{bmatrix} 
\end{CD},  
\]
which is now easily checked as before. 

A similar computation works for the $\cB$-$\cB$ linearity. 
For example, the commutativity of 
\[
\begin{CD}
G(W)\otimes \bB @>>> G(W)\otimes X^*\otimes X\\ 
@VVV @VVV\\ 
\bB\otimes W @>>> X^*\otimes X\otimes W
\end{CD} 
\]
is reduced to that of 
\[
\begin{CD}
W\otimes V\otimes V^* @>>> 
\bigoplus_U U\otimes 
\begin{bmatrix}
U\\ W\,V 
\end{bmatrix} 
\otimes V^*\\ 
@VVV @VVV\\ 
W @<<< \bigoplus_U U\otimes U^*\otimes W, 
\end{CD}
\]
which holds if we define the morphism $\bB \to X^*\otimes X$ 
without weights. 
\end{proof}

\begin{Lemma}
The morphisms $\e: X\otimes X^* \to \bA$ and 
$\d: \bB \to X^*\otimes X$ are supported by 
$e_\cB$ and $e_\cA$ respectively. 
\end{Lemma}

\begin{proof}
We shall check $\e\circ e_\cB = \e$. 
By the commutativity of left and right actions on $X$, 
we see that the composition 
$\sum_k \e\circ \pi({\widehat w}_{kk})$ is given by 
\begin{align*}
X\otimes X^* 
&\to \bigoplus_V X\otimes V\otimes V^*\otimes X^*\\ 
&\to \bigoplus_V F(V)\otimes X\otimes X^*\otimes V^*\\ 
&\to 
\bigoplus_V F(V)\otimes X\otimes G(W)^*\otimes G(W)\otimes X^*\otimes V^*\\ 
&\to \bigoplus_V F(V)\otimes W^*\otimes X\otimes X^*\otimes W|otimes V^*\\ 
&\stackrel{\e_X}{\longrightarrow}
\bigoplus_V F(V)\otimes W^*\otimes W\otimes V^*\\ 
&\stackrel{\sum \pi({\widehat w}_{kk})}{\longrightarrow}
\bigoplus_V F(V)\otimes V^*. 
\end{align*}
From the definition of 
$G(W)\otimes X^* \to X^*\otimes G(W)$, 
the morphism 
\[
\begin{CD}
X\otimes X^* @>{1\otimes \d_{G(W)}\otimes 1}>> 
X\otimes G(W)^*\otimes G(W)\otimes X^* 
@>>> W^*\otimes X\otimes X^*\otimes W 
@>{\e_{W^*X}}>> I 
\end{CD} 
\]
is equal to $d(W) \e_X$. 
Since 
$\sum_k \pi({\widetilde w}_{kk}) = d(W) (1\otimes \e_{W^*}\otimes 1)$, 
we obtain the relation 
\[
\sum_k \e\circ \pi({\widehat w}_{kk}) = d(W)^* \e 
\]
and hence $\e\circ e_\cB = \e$ by taking the summation over the set 
$\{ [W]\}$. 
\end{proof}

We shall now compute 
\[
\begin{CD}
X @>{\o_B^{-1}}>> X\otimes_\cB \bB 
@>{1\otimes \d}>> X\otimes_\cB X^*\otimes_\cA X 
@>{\e\otimes 1}>> \bA\otimes_\cA X 
@>{\o_A}>> X. 
\end{CD} 
\]
As $\e$, $\d$ and $(\lambda,\rho)$ are supported by $e_\cA$ or 
$e_\cB$, the problem is equivalent to consider 
\[
\begin{CD}
X @>{\o_B^{-1}}>> X\otimes \bB 
@>{1\otimes \d}>> X\otimes X^*\otimes X 
@>{\e\otimes 1}>> \bA\otimes X 
@>{\o_A}>> X. 
\end{CD} 
\]
From definition, the composition 
$X \to X\otimes \bB \to X\otimes X^*\otimes X$ is given by 
\begin{align*}
X &\stackrel{\text{weight}}{\longrightarrow} 
\bigoplus_W W^*\otimes W\otimes X 
\to \bigoplus_W X\otimes W^*\otimes G(W) 
\to \bigoplus_W X\otimes W^*\otimes X^*\otimes X\otimes G(W)\\ 
&\longrightarrow \bigoplus_W X\otimes W^*\otimes X^*\otimes W\otimes X 
\to X\otimes X^*\otimes X,  
\end{align*} 
where $\text{weight} = d(W)\o_B^{-1} \dim(B)^{-1}$. By Lemma~6.1, 
this is equivalent to 
\begin{align*}
X &\stackrel{\text{weight}}{\longrightarrow} 
\bigoplus_W W^*\otimes W\otimes X 
\to \bigoplus_W W^*\otimes X\otimes G(W) 
\to \bigoplus_W W^*\otimes X\otimes G(W)\otimes X^*\otimes X\\ 
&\longrightarrow \bigoplus_W W^*\otimes X\otimes X^*\otimes W\otimes X 
\to X\otimes X^*\otimes X.   
\end{align*} 

Similarly, the composition 
$X\otimes X^*\otimes X \to \bA\otimes X \to X$ is given by 
\begin{align*}
XX^*X &\stackrel{\text{weight}}{\longrightarrow} 
\bigoplus_V XV^*VX^*X 
\to \bigoplus_V XV^*X^*F(V)X\\ 
&\longrightarrow \bigoplus_V XV^*X^*XV 
\to XX^*X \stackrel{\e\otimes 1}{\longrightarrow} X
\end{align*}
with $\text{weight} = d(V) \o_A$. 

Note here that by the commutativity $\cT\otimes \cV = \cV\otimes \cT$, 
the position of vector spaces such as $V$ can be freely moved left and right, 
which is pictorially reflected in crossing lines (cf.~Fig). 

Now, combining these two expressions and then applying the definition 
of the trivialization isomorphisms 
$G(W)X^* \to X^*W$, $VX^* \to X^*F(V)$, we have the morphism 
\begin{align*}
X &\to W^*WX \to W^*XG(W) \to W^*F(V)F(V)^*XG(W)\\ 
&\to F(V)W^*XV^*G(W) \to F(V)XG(W^*)V^*G(W) 
\to XVV^*G(W)^*G(W) \to X, 
\end{align*}
which is summed over $[V]$ and $[W]$ with the weight 
$d(V)d(W)\o_A/\o_B \dim(A)$ multiplied (Fig.~\ref{bigfig}). 
By the commutativity of left and right actions, 
we can replace the part $F(V)^*WX \to XV^*G(W)$ with 
\[
WF(V^*)X \to WXV^*\to XG(W)V^* 
\]
to get the expression (Fig.~\ref{longfig})
\begin{align*}
X &\to F(V)W^*WF(V)^*X \to F(V)W^*WXV^* 
\to F(V)W^*XG(W)V^*\\ 
&\to F(V)XG(W^*)G(W)V^* 
\to XVG(W)^*G(W)V^* \to X. 
\end{align*} 

By the associativity of the right action on $X$, the last local morphism 
is reduced to 
\[
X \to F(V)F(V)^*X \to F(V)XV^* \to XVV^* \to X 
\]
multiplied by $d(W)$, which is further reduced to 
$d(V)d(W) 1_X$ by the associativity of the left action on $X$. 

In total, the morphism $X \to XX^*X \to X$ in question amounts to the 
scalar multiple of $1_X$ by 
\[
\sum_{V,W} \frac{d(V)^2d(W)^2}{\dim B} \frac{\o_A}{\o_B} 
= \dim(A) \frac{\o_A}{\o_B}. 
\]

\begin{figure}[htbp]
\vspace*{0.5cm}
\hspace*{-1.5cm}
\input{bigfig.tpc}
\vspace*{0.5cm}
\caption[]{}
\label{bigfig}
\end{figure}

\begin{figure}[htbp]
\vspace*{0.5cm}
\hspace*{-1.5cm}
\input{longfig.tpc}
\vspace*{0.5cm}
\caption[]{}
\label{longfig}
\end{figure}

Similarly, we compute the composition 
\[
\begin{CD}
X^* @>{\o_B^{-1}}>> \bB\otimes X^* 
@>{\d\otimes 1}>> X^*\otimes X\otimes X^* 
@>{1\otimes \e}>> X^*\otimes \bA 
@>{\o_A}>> X^*
\end{CD} 
\]
and see that its is a scalar multiple of $1_{X^*}$ by the same scalar. 

\begin{Proposition}
Let $\cT$ be a rigid semisimple tensor category. 
Then the bicategory $\cM(\cT)$ is rigid as well. 
More precisely, if the unit constraints are specified by 
a function $\{ \o_A \}_A$ indexed by finite-dimensional Hopf algebras 
realized inside $\cT$, then a rigidity pair for an $\cA$-$\cB$ module $X$ 
is given by $(\e,c \d)$ with $c = \dim(A)\o_A/\o_B$, where $\e$ and $\d$ 
are defined above. 
\end{Proposition}

If the tensor category $\cT$ is furnished with a Frobenius duality 
$\{ \e_X: X\otimes X^* \to I \}$ 
(the conjugation being assumed to be strict, which particularly means 
$(X\otimes Y)^* = Y^*\otimes X^*$ and $X^{**} = X$), it is natural to 
use the following normalization for the trivializing isomorphisms of 
the unit object $\bA$: Let the action be renormalized by the gauge 
$\h = \{ \sqrt{d(V)} 1_{V^*} \}$. The morphisms 
$\e: X\otimes X^* \to \bA$ and $\d: \bB \to X^*\otimes X$ are changed into 
the ones associated to the pairing 
\[
V\otimes V^* \ni v\otimes v^* \mapsto \sqrt{d(V)} 
\langle v, v^* \rangle
\]
or its dualized copairing 
\[
\sqrt{d(V)} \sum_i v_i\otimes v_i^* \in V\otimes V^*. 
\]

\begin{Proposition}
Suppose that the semisimple tensor category $\cT$ is furnished 
with a Frobenius duality 
$\{ \e_X \}$ and let the unit constraint $\bA\otimes X \to X$ be 
renormalized by the factor $\o_A = |A|^{-1/2}$ for each $A$ 
with $|A| = \dim A$. 
Then the remormalized family $\{ |A|^{-1/4} |B|^{-1/4} \e \}$ 
gives a Frobenius duality in the bicategory $\cM(\cT)$. 
\end{Proposition} 

\begin{Corollary}[Dimension Formula] 
For an $\cA$-$\cB$ module ${}_\cA X_\cB$, its dimension is calculated by 
\[
\dim({}_\cA X_\cB) = \frac{\dim(X)}{|A|^{1/2} |B|^{1/2}}. 
\] 
Here $\dim(X)$ denotes the dimension of $X$ as an object of $\cT$. 
\end{Corollary}

 
\section{Duality for Orbifolds on Tensor Categories} 

Let $H$ be an off-diagonal object in a rigid semisimple bicategory 
and assume that $H$ satisfies the condition 
\[
H\otimes H^*\otimes H \cong H\oplus \dots \oplus H. 
\] 

Given an object $H$ of this type, we can associate a Hopf algebra 
$B$ so that its Tannaka dual $\cB$ is isomorphic to 
the tensor category generated by $H^*\otimes H$ 
(\cite[Appendix~C]{GSTC}). 
More explicitly, for each object $X$ in $(H^*\otimes H)^n$ 
with $n$ a positive integer, we can construct the monoidal functor 
$X \mapsto E(X)$, where $E(X)$ denotes a finite-dimensional vector space 
defined by  
\[
E(X) = \Hom(H,H\otimes X)
\]
and the multiplicativity isomorphism $E(X)\otimes E(Y) \to E(X\otimes Y)$ 
is given by 
\[
E(X)\otimes E(Y) \ni x\otimes y \mapsto 
(x\otimes 1_Y)y \in E(X\otimes Y). 
\]

\begin{Example} 
Consider the Tannaka dual $\cA$ of a finite-dimensional Hopf algebra $A$ 
realized in a semisimple tensor category $\cT$ and let 
$\bA$ be the associated unit object for $\cA$-$\cA$ modules. 

Then the right $\cA$-module $H = \bA_\cA$ satisfies the above condition. 
In fact, we have 
\[
H\otimes_\cA H^* = \bA = \bigoplus_V V^*\otimes F(V) 
\]
and therefore 
\[
\bigoplus_V V^*\otimes F(V)\otimes \bA_\cA \cong 
\bigoplus_V V^*\otimes \bA_\cA \otimes V 
= \bigoplus_V V^*\otimes V\otimes \bA_\cA  
\] 
is isomorphic to a direct sum of $H$'s. 

Moreover we can identify the associated Hopf algebra with $A$: 
Given an object $V$ in $\cA$, the vector space 
$E(F(V)) = \Hom(H, F(V)\otimes H)$ is naturally isomorphic to $V$ 
by the trivialization isomorphism 
$F(V)\otimes H \cong H\otimes V$ and the simplicity of 
$H_\cA$. Moreover, we have the commutative diagram 
\[
\begin{CD}
V\otimes W @= V\otimes W\\ 
@VVV @VVV\\ 
E(F(V))\otimes E(F(W)) @>>> E(F(V\otimes W))
\end{CD} 
\]
and the monoidal functor $E$ is equivalent to 
the identity functor in $\cA$. Thus the associated Hopf algebra 
is naturally isomorphic to $A$, whereas the object $H^*\otimes H$ 
generates the tensor category isomorphic to the Tannaka dual of 
the dual Hopf algebra $B = A^*$. 
\end{Example}

\begin{Proposition}
The construction of Hopf algebras 
from objects of absorbing property is universal, i.e., 
any finite-dimensional semisimple Hopf algebra arises this way. 
\end{Proposition}

Returning to the initial case of this section, the obvious identification 
\[
H\otimes X \to E(X)\otimes H
\]
can be interpreted as giving a right action of $\cB$ on $H$. 

Consider the composite isomorphism 
\[
H^*\otimes H \to \bigoplus_X X\otimes \Hom(X,H^*\otimes H) 
\to \bigoplus_X X\otimes E(X^*) = \bB. 
\]
We shall show that this isomorphism is $\cB$-$\cB$ linear, 
i.e., the commutativity of 
\[
\begin{CD}
X\otimes H^*\otimes H\otimes Y @>>> X\otimes \bB\otimes Y\\ 
H^*\otimes E(X)\otimes E(Y)\otimes H @>>> 
E(Y)\otimes \bB\otimes E(X)  
\end{CD} 
\]
or equivalently, by applying the functor $\Hom(Z,\cdot)$ with 
$Z$ a simple object, we have the commutative diagram of vector spaces. 
For simplicity, letting $X = I$ (the letter $X$ will be used as 
a dummy index), the relevant isomorphisms are given by 
\[
\begin{CD}
  \begin{bmatrix}
    Z\\ H^*HY
  \end{bmatrix} 
@>>> 
\bigoplus_X 
\begin{bmatrix}
  X\\ H^*H
\end{bmatrix} 
\otimes 
\begin{bmatrix}
  Z\\ X\,Y 
\end{bmatrix} 
@>>> \bigoplus_X 
\begin{bmatrix}
  H\\ H\,X^* 
\end{bmatrix} 
\otimes 
\begin{bmatrix}
  X^*\\ Y\,Z^* 
\end{bmatrix}\\ 
@VVV @. @VVV\\ 
\begin{bmatrix}
  H\\ H\,Y 
\end{bmatrix} 
\otimes 
\begin{bmatrix}
  Z\\ H^*\,H 
\end{bmatrix} 
@>>> 
\begin{bmatrix}
  H\\ H\,Y 
\end{bmatrix} 
\otimes 
\begin{bmatrix}
  H\\ H\,Z^* 
\end{bmatrix} 
@>>> 
\begin{bmatrix}
  H\\ HYZ^* 
\end{bmatrix}. 
\end{CD} 
\]

To check the commutativity, let us start with a vector 
$x\otimes T \in 
\begin{bmatrix}
H\\ H\,X^* 
\end{bmatrix} 
\otimes 
\begin{bmatrix}
X^*\\ Y\,Z^* 
\end{bmatrix}$. 
The upper horizontal line is then described as 
\[
({\widetilde x}\otimes 1){\widetilde T} 
\mapsto {\widetilde x}\otimes {\widetilde T} 
\mapsto x\otimes T, 
\]
while the right and the left vertical lines are presented by 
$x\otimes T \mapsto (1\otimes T)x$ and 
\[
({\widetilde x}\otimes 1){\widetilde T} 
\mapsto 
\sum_{j,k} 
\langle z_k^*(1\otimes y_j^*), 
({\widetilde x}\otimes 1){\widetilde T} \rangle 
y_j\otimes z_k 
\]
with $\{ y_j, y_j^*\}$ and $\{ z_k, z_k^*\}$ in the duality relation 
($y_j^*y_j = 1_H$ and $z_k^*z_k = 1_Z$ particularly). 
Finally the bottom line is given by 
\[
\sum_{jk} c_{jk} y_j\otimes z_k 
\mapsto 
\sum_{j,k} c_{jk} y_j\otimes {\widetilde z}_k 
\mapsto 
\sum_{j,k} c_{jk} (y_j\otimes 1) {\widetilde z}_k. 
\]

To identify the last summation with $(1\otimes T)x$, 
we rewrite $c_{jk}$ as follows: 
\begin{align*}
d(Z)c_{jk} &= 
\e_Z(z_k^*\otimes 1)(1\otimes y_j^*\otimes 1) 
({\widetilde x}\otimes 1)({\widetilde T}\otimes 1)\d_{Z^*}\\ 
&= \e_Z(z_k^*\otimes 1)(1\otimes y_j^*\otimes 1) 
(1\otimes \e_{X^*}\otimes 1)(1\otimes x\otimes {\widetilde T}\otimes 1) 
(\d_H\otimes \d_{Z^*})\\ 
&= \e_{H^*}(1\otimes \widetilde{z_k^*}) 
(1\otimes y_j^*\otimes 1)(1\otimes \e_{X^*}\otimes 1) 
(1\otimes x\otimes {\widetilde T}\otimes 1) 
(1\otimes \d_{Z^*})\d_H,  
\end{align*}
which yields the relation 
\[
\frac{d(Z)}{d(H)} c_{jk} 1_H 
= \widetilde{z_k^*}(y_j^*\otimes 1)(1\otimes \e_{X^*}\otimes 1) 
(x\otimes {\widetilde T}\otimes 1)(1_H\otimes \d_{Z^*}). 
\]
Now this formula is used to get 
\[
\sum_k c_{jk} \widetilde{z_k} = 
\sum_k c_{jk} \widetilde{z_k} 1_H = 
\sum_k \frac{d(H)}{d(Z)} 
\widetilde{z_k} \widetilde{z_k^*} (y_j^*\otimes 1) 
(1\otimes \e_{X^*}\otimes 1) (x\otimes {\widetilde T}\otimes 1) 
(1_H\otimes \d_{Z^*}). 
\]

From the relation 
\[
\langle \widetilde{z_k^*} \widetilde{z_k} \rangle 
= \e_Z(z_k^*\otimes 1)(z_k\otimes 1)\d_{Z^*} = d(Z), 
\]
we see that 
$\widetilde{z_k^*} \widetilde{z_k} = d(Z)/d(H) 1_H$ 
and hence 
\[
(\widetilde{z_k})^* = \frac{d(H)}{d(Z)} \widetilde{z_k^*}. 
\]
Feeding this back into the above summation, we have 
\[
\sum_k c_{jk} \widetilde{z_k} = 
(y_j^*\otimes 1)(1\otimes \e_{X^*}\otimes 1) 
(x\otimes {\widetilde T}\otimes 1) (1_H\otimes \d_{Z^*}) 
\]
and then 
\begin{align*}
\sum_{j,k} c_{jk} (y_j\otimes 1) \widetilde{z_k} 
&= \sum_j (y_jy_j^*\otimes 1) (1\otimes \e_{X^*}\otimes 1) 
(x\otimes {\widetilde T}\otimes 1) (1_H\otimes \d_{Z^*})\\ 
&= (1\otimes \e_{X^*}\otimes 1) 
(x\otimes {\widetilde T}\otimes 1) (1_H\otimes \d_{Z^*})\\ 
&= (1\otimes T)x. 
\end{align*}

\begin{Lemma}
We have 
\[
{}_\cB H^*\otimes H_\cB \cong {}_\cB \bB_\cB, 
\qquad 
H\otimes_\cB H^* \cong I. 
\]
\end{Lemma}

\begin{proof}
We have just checked the former relation. By Frobenius reciprocity, 
this implies 
\[
\dim \End(H\otimes_\cB H^*) = \dim \End({}_\cB H^*\otimes H_\cB) = 1 
\] 
and hence $H\otimes_\cB H^* = I$ by semisimplicity. 
\end{proof}

Since bimodules with the similar property are referred to 
as imprimitivity bimodules in connection with Mackey's imprimitivity 
theorem on induced representations, we call an object $M$ 
in a rigid bicategory 
an {\bf imprimitivity object} if both of $M\otimes M^*$ and $M^*\otimes M$ 
are isomorphic to unit objects. 
In a tensor category, this is nothing but saying that $M$ is 
an invertible object. 

The following observation, though obvious, is the essence of duality 
for orbifold constructions. 

\begin{Lemma}
Let 
\[
\begin{pmatrix}
\cT & \cM\\ 
\cM^* & \cS
\end{pmatrix} 
\] 
be a rigid semisimple bicategory and $M$ be an imprimitivity object 
in $\cM$. 

Then two tensor categories $\cS$ and $\cT$ are isomorphic. 
More precisely, 
\[
X \mapsto M\otimes X\otimes M^*, 
\qquad 
Y \mapsto M^*\otimes Y\otimes M 
\]
gives the monoidal equivalence between $\cS$ and $\cT$. 
\end{Lemma}

Given a monoidal imbedding $F: \cA \to \cT$ of the Tannaka dual 
$\cA$ of a finite-dimensional semisimple Hopf algebra $A$ 
into a rigid semisimple tensor category $\cT$, let 
$H = \bA_\cA$ be an off-diagonal object in the bicategory 
\[
\begin{pmatrix}
\cT & \cM_\cA\\ 
{}_\cA \cM & {}_\cA \cM_\cA. 
\end{pmatrix}
\]
Here $\cM_\cA$ denotes the category of right $\cA$-modules in $\cT$ 
and similarly for others. 

Then $H$ meets the absorbing property and the tensor subcategory 
of ${}_\cA \cM_\cA$ 
generated by $H^*\otimes H = {}_\cA \bA\otimes \bA_\cA$ is 
isomorphic to the Tannaka dual $\cB$ of the dual Hopf algebra of $A$.  
Let $G: \cB \to {}_\cA \cM_\cA$ be the accompanied monoidal imbedding. 
Recall here that the Tannaka dual $\cA$ of $A$ is the one associated 
to $H\otimes H^*$ as seen in the above example. 

Thus we can talk about $\cB$-modules in $\cM$: 
Let $\cM_\cB$ (resp.~${}_\cB \cM$) be the category of 
right (resp.~left) $\cB$-modules in 
$\cM_\cA$ (resp.~${}_\cA \cM$) and ${}_\cB \cM_\cB$ be the 
category of $\cB$-$\cB$ bimodules in ${}_\cA \cM_\cA$. 
Then these, together with the starting tensor category $\cT$, 
form a bicategory 
\[
\begin{pmatrix}
\cT & \cM_\cB\\ 
{}_\cB \cM & {}_\cB \cM_\cB. 
\end{pmatrix} 
\]
Thanks to the previous discussions, 
the object $H = \bA_\cA$ in $\cM_\cA$ admits a structure of 
right $\cB$-module, which gives rise to an imprimitivity object 
$M_\cB$ in $\cM_\cB$. Then the above lemma shows that the tensor 
category ${}_\cB \cM_\cB$ is isomorphic to the original tensor 
category. 

To get the meaning of this, we first introduce the notation 
$\cT\rtimes_F\cA$ for the tensor category ${}_\cA\cM_\cA$, which is 
interpreted as the crossed product of $\cT$ by $F$. 
Then the monoidal imbedding $G: \cB \to \cT\rtimes_F\cA$ describes 
the dual symmetry in $\cT\rtimes_F\cA$ and we can construct 
the second crossed product $(\cT\rtimes_F\cA)\rtimes_G\cB$. 

\begin{Theorem}
With the notation described above, we have the duality for crossed products: 
the second crossed product 
tensor category $(\cT\rtimes_F\cA)\rtimes_G\cB$ is canonically 
isomorphic to the original tensor category $\cT$. 
\end{Theorem}


 

\begin{thebibliography}{42}
\bibitem{CP}
V.~Chari and A.~Pressley, 
{\it A Guide to Quantum Groups}, 
Cambridge Univ.~Press, 
1995.
\bibitem{Da}
A.A.~Davydov, 
Monoidal categories,  
{\it J.~Math.~Sci.~(New York)}, 
88(1998), 457--519.
\bibitem{DZ}
P.~Di Francesco and J.-B.~Zuber, 
$SU(N)$ lattice integrable models associated with graphs, 
{\it Nuclear Physics B}, 
338(1990), 602--646. 
\bibitem{DPR} 
S.~Doplicher, C.~Pinzari and J.E.~Roberts, 
An algebraic duality theory for multiplicative unitaries, 
preprint, 2000.  
\bibitem{EK1}
D.~Evans and Y.~Kawahigashi, 
Orbifold subfactors from Hecke algebras
{\it Commun.~Math.~Phys.}, 
165(1994), 445--484.
\bibitem{EK2}
D.~Evans and Y.~Kawahigashi, 
{\it Quantum Symmetries on Operator Algebras}, 
Oxford University Press, 
Oxford, 
1998. 
\bibitem{F}
J.M.G.~Fell, 
An extension of Mackey's method to Banach *-algebraic bundles, 
{\it Memoirs Amer.~Math.~Soc.}, 
90(1969). 
\bibitem{FG}
P.~Fendley and P.~Ginsparg, 
Non-critical orbifolds, 
{\it Nuclear Physics B}, 
324(1989), 549--580. 
\bibitem{FK}
J.~Fr\"ohlich and T.~Kerler, 
{\it Quantum Groups, Quantum Categories and Quantum Field Theory}, 
Lecture Notes in Math.~1542, Springer-Verlag, Berlin, 1993. 
\bibitem{Kas}
C.~Kassel, 
{\it Quantum Groups}, 
Springer-Verlag, 
Berlin-New York, 
1995. 
\bibitem{KaY}
T.~Kajiwara and S.~Yamagami, 
Irreducible bimodules associated with crossed product algebras II, 
{\it Pacific J.~Math.},
171(1995), 209--229.
\bibitem{Ki}
A.~A.~Kirillov, 
{\it Elements of the Theory of Representations}, 
Springer-Verlag, 
Berlin, 
1976.
\bibitem{KMY}
H.~Kosaki, A.~Munemasa and S.~Yamagami,
Fusion algebras associated to finite group actions,
{\it Pacific J. Math.},
117(1997), 269--290.
\bibitem{KoY}
H.~Kosaki and S.~Yamagami,
Irreducible Bimodules Associated with Crossed Product Algebras,
{\it Int. J. Math.},
3(1992), 661--676.
\bibitem{LaRa1}
R.G.~Larson and D.E.~Radford, 
Semisimple cosemisimple Hopf algebras, 
{\it Amer.~J.~Math.}, 
109(1987), 187--195.  
\bibitem{LaRa2}
R.G.~Larson and D.E.~Radford, 
Finite dimensional cosemisimple Hopf algebras in characteristic $0$ 
are semisimple, 
{\it J.~Algebra}, 
117(1988), 267--289. 
\bibitem{L}
R.~Longo, 
A duality for Hopf algebras and for subfactors. I, 
{\it Commun.~Math.~Phys.}, 
159(1994), 133--150.
\bibitem{Mac}
S.~MacLane, 
{\it Categories for the Working Mathematician}, 
Springer-Verlag, 
Berlin-New York, 
1971. 
\bibitem{Maj}
S.~Majid, 
Reconstruction theorems and rational conformal field theories, 
{\it Internat.~J.~Modern Phys.~A}, 
6(1991), 4359--4374. 
\bibitem{Mal}
G.~Maltsiniotis, 
Traces dans les cat\'egories monoidales, dualit\'e et cat\'egories monoidales fibries, 
{\it Cahiers Topologie G\'eom. Diff\'erentielle Cat\'eg.},
36(1995), 195--288.
\bibitem{R} 
M.A.~Rieffel, 
Induced representations of C*-algebras, 
{\it Adv.~Math.}, 
13(1974), 176--257. 
\bibitem{SS} 
S.~Shnider and S.~Sternberg, 
{\it QuantumGroups}, 
International Press, 
Boston, 
1993. 
\bibitem{Sz} 
W.~Szymanski, 
Finite index subfactors and Hopf algebra crossed products, 
{\it Proc.~Amer.~Math.~Soc.},
120(1994), 519--528. 
\bibitem{T}
M.~Takeuchi, 
Matched pairs of groups and bismash products of Hopf algebras, 
{\it Commun.~Algebra}, 
9(1981), 841--882. 
\bibitem{Tam} 
D.~Tambara, 
A duality for modules over monoidal categories of representations of 
semisimple Hopf algebras, 
preprint, 
1998. 
\bibitem{TY}
D.~Tambara and S.~Yamagami, 
Tensor categories with fusion rules of self-duality for 
finite abelian groups, 
{\it J.~Algebra}, 
209(1998), 692--707
\bibitem{U}
K.-H.~Ulbrich, 
On Hopf algebras and rigid monoidal categories, 
{\it Israel J.~Math.}, 
72(1990), 252--256. 
\bibitem{W}
S.L.~Woronowicz, 
Tannaka-Krein duality for compact matrix pseudogroups. 
Twisted $SU(N)$ groups, 
{\it Invent.~Math.}, 
93(1988), 35--76. 
\bibitem{VB}
S.~Yamagami, 
Vector bundles and bimodules,
{\it Proceedings in Quantum and Non-Commutative Analysis},
Kluwer Academic,
1993, 321--329. 
\bibitem{URQG}
\underline{\phantom{S.~Yamagami}}, 
On unitary representation theories of compact quantum groups, 
{\it Commun.~Math.~Phys.}, 
167(1995), 509--529. 
\bibitem{FRTC}
\underline{\phantom{S.~Yamagami}}, 
Frobenius reciprocity in tensor categories, 
{\it Math.~Scand.}, 
to appear.
\bibitem{GSTC}
\underline{\phantom{S.~Yamagami}}, 
Group symmetry in tensor categories and duality for orbifolds, 
preprint, 2000.
\bibitem{FPB}
\underline{\phantom{S.~Yamagami}}, 
C*-tensor categories and free product bimodules, 
preprint, 
1999. 
\bibitem{FDTC}
\underline{\phantom{S.~Yamagami}}, 
Frobenius duality in C*-tensor categories, 
preprint, 
2000. 
\end{thebibliography}
\end{document}